\documentclass{gtart_h}  

\def\ifplaintex{\expandafter\ifx\csname documentclass\endcsname\relax}

\def\ifplaintex{\expandafter\ifx\csname documentclass\endcsname\relax}

\ifplaintex 
\else
\fi

\expandafter\ifx\csname epsfbox\endcsname\relax\input epsf\fi

\def\gt{{\mathsurround=0pt\it $\cal G\mskip-2mu$eometry \&\ 
$\cal T\!\!$opology}}        %  journal title in recommended style

\def\gtp{{\mathsurround=0pt\it $\cal G\mskip-2mu$eometry \&\ 
$\cal T\!\!$opology $\cal P\!$ublications}}  % GT publications

\def\lognumber#1{\def\thelognumber{#1}}
\def\volumenumber#1{\def\thevolumenumber{#1}}
\def\papernumber#1{\def\thepapernumber{#1}}
\def\volumeyear#1{\def\thevolumeyear{#1}}

\def\pagenumbers#1#2{\def\startpage{#1}\def\finishpage{#2}}
\def\published#1{\def\publishdate{#1}}
\def\proposed#1{\def\theproposer{#1}}
\def\seconded#1{\def\theseconders{#1}}
\def\received#1{\def\receiveddate{#1}}

\def\accepted#1{\def\accepteddate{#1}}

\let\\\par\let\thelognumber\relax
\let\thevolumenumber\relax\let\thepapernumber\relax
\let\thevolumeyear\relax\let\thesamplenumber\relax\let\startpage\relax
\let\finishpage\relax\let\publishdate\relax\let\receiveddate\relax
\let\reviseddate\relax\let\accepteddate\relax\let\theasciititle\relax
\let\theasciiauthors\relax
\let\theasciiabstract\relax
\let\theasciiemail\relax\let\theshortauthors\relax\let\theshorttitle\relax

\long\def\maketitlep{   % start of definition of \maketitlep

\count0=\startpage

\gt\hfill      %   Journal title (top left) 
\hbox to 77pt{\vbox to 0pt{\vglue -15pt\epsfbox{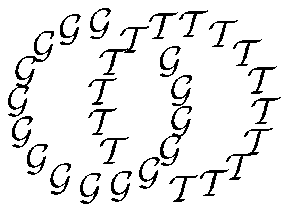}\vss}\hss}
\break
{\small\ifx\thesamplenumber\relax % sample?  
Volume \else Sample
\fi\thevolumenumber\ (\thevolumeyear)
\startpage--\finishpage\nl
Published: \publishdate}
\vglue 0.5truein plus 0.4fil minus 0.1truein

{\parskip=0pt\leftskip 0pt plus 1fil\def\\{\par\smallskip}{\ifplaintex\large
\else\Large\fi\bf\thetitle}\par\medskip}   

\vglue 0pt plus 0.1fil 

{\parskip=0pt\leftskip 0pt plus 1fil\def\\{\par}{\sc\theauthors}
\par\medskip}

\vglue 0pt plus 0.1fil 

{\small\parskip=0pt\let\newline\\
{\leftskip 0pt plus 1fil\def\\{\par}{\sl\theaddress}\par}
\expandafter\ifx\theemail\relax    % email address?
\relax\else\vglue 5pt plus 0.02fil minus 2pt\def\\{\stdspace{\rm 
and}\stdspace} 
\cl{Email:\stdspace\tt\theemail}\fi
\ifx\theurl\relax                  % URL given?
\relax\else\vglue 5pt plus 0.02fil minus 2pt\def\\{\stdspace{\rm 
and}\stdspace}
\cl{URL:\stdspace\tt\theurl}\fi\par}

\vglue 7pt plus 0.3fil minus 3pt

{\bf Abstract}
\vglue 5pt plus 0.1fil minus 2pt

\theabstract

\vglue 7pt plus 0.3fil minus 3pt

{\bf AMS Classification numbers}\quad Primary:\quad \theprimaryclass

Secondary:\quad \thesecondaryclass

\vglue 5pt plus 0.3fil minus 2pt

{\bf Keywords:}\quad \thekeywords

\vglue 10pt plus 0.5fil minus 5pt

{\small  Proposed: \theproposer\hfill Received: \receiveddate\nl
Seconded: \theseconders\hfill 
\ifx\reviseddate\relax                         % paper revised?
Accepted: \accepteddate                        % no
\else
Revised: \reviseddate                          % yes
\fi}
\eject
}       %  end of definition of \maketitlep

\font\phead=cmsl9 scaled 950
\font=cmsl9 scaled 1050
\font\pnum=cmbx10 scaled 913
\font\lnum=cmbx10 
\font\pfoot=cmsl9 scaled 950
\font=cmsl9 scaled 1050
\ifplaintex
\headline{\vbox to 0pt{\vskip -4.5mm\line{\small\phead\ifnum
\count0=\startpage ISSN 1364-0380 (on line)
1465-3060 (printed) \hfill {\pnum\folio}\else\ifodd\count0\def\\{ }% 
\ifx\theshorttitle\relax\thetitle\else\theshorttitle\fi\hfill{\pnum\folio}
\else\def\\{ and }{\pnum\folio}\hfill\ifx\theshortauthors\relax\theauthors
\else\theshortauthors\fi\fi\fi}\vss}}
\footline{\vbox to 0pt{\vglue 0mm\line{\small\pfoot\ifnum\count0=\startpage
\copyright\ \gtp\hfill\else
\gt, Volume \thevolumenumber\ (\thevolumeyear)\hfill\fi}\vss
}}
\else
\makeatletter
\def\@oddhead{{\small\lhead\ifnum\count0=\startpage ISSN 1364-0380 (on line)
1465-3060 (printed) \hfill {\lnum\number\count0}\else\ifodd\count0
\def\\{ }\ifx\theshorttitle\relax \thetitle \else\theshorttitle\fi\hfill
{\lnum\number\count0}\else\def\\{ and }{\lnum\number\count0}
\hfill\ifx\theshortauthors\relax 
\theauthors\else\theshortauthors\fi\fi\fi}}\def\@evenhead{\@oddhead}
\def\@oddfoot{\small\lfoot\ifnum\count0=\startpage\copyright\ \gtp\hfill\else
\gt, Volume \thevolumenumber\ (\thevolumeyear)\hfill\fi}
\def\@evenfoot{\@oddfoot}
\makeatother
\fi

\newwrite\gtoutfile
\long\gdef\makeheadfile{  %%% start of definition of \makeheadfile
{\def\\{, }\def\s{ }
\immediate\openout\gtoutfile head.xxx
\immediate\write\gtoutfile{Proxy-for: \ifx\theasciiauthors\relax
\theauthors\else\theasciiauthors\fi\s<\ifx\theasciiemail\relax\theemail\else\theasciiemail\fi>}
\immediate\write\gtoutfile{\noexpand\\}
\immediate\write\gtoutfile{Authors: \ifx\theasciiauthors\relax
\theauthors\else\theasciiauthors\fi}
{\def\\{ }\immediate\write\gtoutfile{Title: \ifx\theasciititle\relax
\thetitle\else\theasciititle\fi}}
\immediate\write\gtoutfile{Subj-class: GT or SG or MG etc}
\immediate\write\gtoutfile{MSC-class: \theprimaryclass\ifx\thesecondaryclass\relax\else, \thesecondaryclass\fi}
\immediate\write\gtoutfile{Journal-ref: Geom. Topol. \thevolumenumber
(\thevolumeyear) \startpage-\finishpage}
\immediate\write\gtoutfile{Comments: Published by Geometry and Topology at}
\immediate\write\gtoutfile{\s\s http://www.maths.warwick.ac.uk/gt/GTVol\thevolumenumber/paper\thepapernumber.abs.html}
\immediate\write\gtoutfile{\noexpand\\}
\immediate\write\gtoutfile{}
\ifx\theasciiabstract\relax
\immediate\write\gtoutfile{\theabstract}\else
\immediate\write\gtoutfile{\theasciiabstract}\fi
\immediate\write\gtoutfile{}
\immediate\write\gtoutfile{\noexpand\\}
\immediate\write\gtoutfile{}
\immediate\closeout\gtoutfile}}  %%% end of definition of \makeheadfile

\def\maketitlepage{\maketitlep\makeheadfile}
\let\maketitle\maketitlepage

\lognumber{478}
\received{16 July 2004}
\volumenumber{9}\papernumber{8}\volumeyear{2005}
\pagenumbers{247}{297}   
\published{26 January 2005}
\accepted{13 January 2005}
\proposed{Yasha Eliashberg}
\seconded{Robion Kirby, Joan Birman}

\usepackage{amsmath,amssymb}
\allowdisplaybreaks
\input xy
\xyoption{all}

\newtheorem{proposition}{Proposition}
\newtheorem{theorem}[proposition]{Theorem}
\newtheorem{lemma}[proposition]{Lemma}
\newtheorem{corollary}[proposition]{Corollary}

\theoremstyle{definition}
\newtheorem{definition}[proposition]{Definition}

\numberwithin{proposition}{section}

\def\Z{\mathbb{Z}}
\def\R{\mathbb{R}}
\def\Q{\mathbb{Q}}
\def\A{\mathcal{A}}

\def\M{\mathcal{M}}
\def\Lag{\mathcal{L}}
\def\d{\partial}
\def\alg#1{\Z\langle{#1}\rangle}
\def\Aut{\operatorname{Aut}}
\def\dott{*}
\def\disc{D}
\def\hom#1{\phi_{#1}}
\def\varhom#1{\tilde{\phi}_{#1}}
\def\Phil#1{\Phi^L_{#1}}
\def\Phir#1{\Phi^R_{#1}}
\def\varPhil#1{\tilde{\Phi}^L_{#1}}
\def\varPhir#1{\tilde{\Phi}^R_{#1}}
\def\varHC{\widetilde{HC}}
\def\HCknot{HC_0^{\operatorname{knot}}}
\def\Burau#1{{\operatorname{Bur}_{#1}}}
\def\red#1{{\operatorname{Bur}^0_{#1}}}
\def\Bur#1{{\operatorname{Bur}_{#1}}}
\def\im{\operatorname{im}}
\def\ab{\operatorname{ab}}
\def\Aug{\operatorname{Aug}}
\def\Augab{\operatorname{Aug}^{\operatorname{ab}}}
\def\ext{\operatorname{ext}}
\def\lin{\operatorname{lin}}
\def\Sym{\operatorname{Sym}}
\def\I{\mathcal{I}}
\def\braid{\operatorname{braid}}
\def\transpose{T}
\def\knot{\operatorname{knot}}

\def\heighten{\vrule width 0pt depth 3pt height 12pt}

\begin{document}

\title{Knot and braid invariants from contact homology I}
\author{Lenhard Ng}
\address{Department of Mathematics, Stanford
University\\Stanford, CA 94305, USA}
\email{lng@math.stanford.edu}
\urladdr{http://alum.mit.edu/www/ng}

\begin{abstract}
We introduce topological invariants of knots and braid conjugacy
classes, in the form of differential graded algebras, and present an
explicit combinatorial formulation for these invariants. The
algebras conjecturally give the relative contact homology of certain
Legendrian tori in five-dimensional contact manifolds. We present
several computations and derive a relation between the knot
invariant and the determinant.
\end{abstract}


\maketitlepage

%*********************************************************************
%*********************************************************************
\section{Introduction}

The study of invariants of Legendrian submanifolds in contact
manifolds is currently a very active field of research. Holomorphic
curve techniques stemming from the introduction of Symplectic Field
Theory \cite{EGH} have inspired a great deal of work on Legendrian
isotopy invariants; papers in the subject include
\cite{Che,EES,Eli,EFM,ENS,Ng,Sab}, among others. Most of these
papers focus on Legendrian knots, and all deduce results in
symplectic and contact geometry. Legendrian knots, however, form a
restrictive class within the space of all knots, and Legendrian
isotopy only makes sense in the presence of an ambient contact
structure.

In this manuscript, we use ideas from Symplectic Field Theory to
obtain new invariants of topological knots and braids. We note that
holomorphic curves have been previously used to define knot
invariants, through the work of Ozsv\'ath and Szab\'o \cite{OSz} on
Heegaard Floer homology. Their approach, however, is very different
from the one leading to the invariants in our paper, which is due to
Eliashberg and can be briefly summarized as follows. A knot or braid
yields in a natural way a Legendrian torus in a certain
$5$-dimensional contact manifold, and topological isotopy of the
knot or braid results in Legendrian isotopy of the torus; thus,
Legendrian invariants of the torus, and in particular relative
contact homology, yield topological invariants for the original knot
or braid.

We will discuss this correspondence further in
Section~\ref{sec:contact}, but the analytical underpinnings of the
precise calculation of relative contact homology are still in
progress. Instead, in the rest of the paper, we handle our new
invariants algebraically, with no reference to contact geometry. The
proofs of invariance, and all subsequent calculations, are
completely combinatorial and algebraic.

We now describe our invariants. The starting point is a
representation $\phi$ of the braid group $B_n$ as a group of algebra
automorphisms of a tensor algebra with $n(n-1)$ generators. This
representation, which extends ones considered previously by Magnus
\cite{Mag} and Humphries \cite{Hum}, is nearly faithful, and there
is a slightly larger representation $\phi^{\ext}$ which is faithful.
The map $\phi$ then gives rise to a differential graded algebra
(customarily abbreviated ``DGA'' in the subject) for any braid $B$,
which we call the braid DGA, while $\phi^{\ext}$ yields a larger
algebra called the knot DGA.

There is an equivalence relation on DGAs first introduced by
Chekanov \cite{Che} called stable tame isomorphism, which has been
the subject of much work because of its importance in the theory of
Legendrian invariants. The main results of this paper state that, up
to stable tame isomorphism, the braid DGA of $B$ is an invariant of
the conjugacy class of $B$, while the knot DGA of $B$ is an
invariant of the knot which is the closure of $B$.

In particular, the graded homology of these DGAs, which is unchanged
by stable tame isomorphism, gives an invariant of braid conjugacy
classes and knots; we call this the contact homology of the braid or
knot. Although the homology of a DGA is generally hard to compute,
there is a relatively simple form for the $0$-dimensional part
$HC_0$ of the contact homology, and we can compute $HC_0$ explicitly
for a large number of knots.

We concentrate on the invariant $HC_0$ because it is difficult in
general to distinguish between stable tame isomorphism classes of
DGAs, and because $HC_0$ seems to encode much of the information
from the knot or braid DGA. It is sometimes also hard to tell when
two knots or braids have isomorphic $HC_0$, but we can deduce
numerical invariants from $HC_0$, called augmentation numbers, which
are easy to calculate by computer.

Knot contact homology seems to be an invariant distinct from
previously known knot invariants. On the one hand, it does not
distinguish between mirrors, and a result in \cite{II} states that,
for two-bridge knots, $HC_0$ depends only on the determinant
$|\Delta_K(-1)|$ of the knot; hence one might expect that knot
contact homology is fairly weak as an invariant. On the other hand,
there are pairs of knots which are distinguished by $HC_0$ (via
augmentation numbers) but not by (any given one of) a litany of
familiar knot invariants: Alexander polynomial, Jones polynomial,
HOMFLY polynomial, Kauffman polynomial, signature, Khovanov
invariant, and Ozsv\'ath--Szab\'o invariant.

One classical invariant that is contained in knot contact homology
is the determinant; in fact, the determinant is encoded by a
linearization $HC_*^{\lin}$ of $HC_*$. A further link between
$HC_0(K)$ and the determinant is the existence of a surjection from
$HC_0(K)$ to a ring $\Z[x]/(p(x))$, where $p$ is a certain
polynomial determined by the homology of the double branched cover
over $K$; this map is often an isomorphism.

Although knot contact homology is uniquely defined as a theory over
$\Z_2$, it has different viable lifts to $\Z$, one of which is the
theory described above. However, there is another, inequivalent,
lift to $\Z$, which we call alternate knot contact homology; the
alternate theory lacks some of the properties of the usual theory
but produces a knot invariant which seems to contain more
information.

The paper is organized as follows. In Section~\ref{sec:defs}, we
define the braid and knot DGAs and state the main invariance
results. Section~\ref{sec:contact} gives a highly informal
description of how our invariants relate to contact geometry.
Section~\ref{sec:HC} defines contact homology for braids and knots
and proves invariance of $HC_0$; we subsequently return to the
invariance proof of the full DGAs in Section~\ref{sec:fullproof}.
Section~\ref{sec:properties} presents properties of the invariants,
including augmentation numbers and behavior under mirrors, and
Section~\ref{sec:determinant} relates knot contact homology to the
knot's determinant. Section~\ref{sec:computations} gives some
examples of computations of the invariants. We discuss the alternate
knot invariant in Section~\ref{sec:alternate}.

In the sequel to this paper \cite{II}, we give a geometric
interpretation of $\phi$ via the usual treatment of $B_n$ as a
mapping class group. We then use this to present a new, completely
topological definition for $HC_0$ for braids and knots, in terms of
cords and relations similar to skein relations, and find a relation
to character varieties. This interpretation explains, independent of
contact geometry, ``why'' $HC_0$ is a topological invariant.

\rk{Acknowledgements}
I am grateful to Yasha Eliashberg for suggesting the project and the
means to approach it and offering a great deal of insightful
guidance. I would also like to thank John Etnyre, Michael Hutchings,
Kiran Kedlaya, Tom Mrowka, Josh Sabloff, Ravi Vakil, and Ke Zhu for
useful conversations, and the Institute for Advanced Study, the
American Institute of Mathematics, and Stanford University for their
hospitality. Funding was provided by a Five-Year Fellowship from the
American Institute of Mathematics, and at one point by NSF grant
DMS-9729992.

%*********************************************************************
%*********************************************************************
\section{Definitions and main results}
\label{sec:defs}

%*********************************************************************
\subsection{Differential graded algebras}
\label{ssec:DGA}

In contact geometry, the homology of a certain DGA defines the
relative contact homology of a Legendrian submanifold in a contact
manifold. In all known cases, a Legendrian isotopy invariant is
provided not only by relative contact homology, but also by the
stable tame isomorphism class of the DGA.

We recall the basic definitions, originally due to Chekanov
\cite{Che}. Let $\alg{a_1,\dots,a_n}$ denote the (noncommutative,
unital) tensor algebra generated by $a_1,\dots,a_n$, which is
generated as a $\Z$-module by words in the $a_i$'s, including the
empty word.

\begin{definition}
Equip $\A=\alg{a_1,\dots,a_n}$ with a grading over $\Z$ by assigning
degrees to the generators $a_i$. Then a \textit{differential} on
$\A=\alg{a_1,\dots,a_n}$ is a map $\d:\,\A\to\A$ which lowers degree
by $1$ and satisfies $\d^2=0$ and the Leibniz rule $\d(vw) = (\d v)
w + (-1)^{\deg v} v(\d w)$. The pair $(\A,\d)$ is what we call a
\textit{differential graded algebra}, for the purposes of this
paper.
\end{definition}

To define a particular equivalence relation between DGAs known as stable
tame isomorphism, we need several more notions, beginning with a
refinement of the usual notion of isomorphism.

\begin{definition}
An algebra map $\phi\co(\alg{a_1,\dots,a_n},\d)\to
(\alg{\tilde{a}_1,\dots,\tilde{a}_n},\tilde{\d})$ is an
\textit{elementary isomorphism} if it is a graded chain map and, for
some $i$,
\[
\phi(a_j) = \begin{cases} \pm \tilde{a}_i+v & \textrm{for $j=i$
and some $v\in
\alg{\tilde{a}_1,\dots,\tilde{a}_{i-1},\tilde{a}_{i+1},\dots,\tilde{a}_n}$}
\\
\tilde{a}_j & \textrm{for $j\neq i$}.
\end{cases}
\]
A \textit{tame isomorphism} between DGAs is any composition of
elementary isomorphisms.
\end{definition}
In the above definition, the order of the generators is
immaterial; that is, a map which simply permutes the generators of
the DGA (and treats the grading and differential accordingly) is
considered to be elementary. Note also that elementary and tame
isomorphisms are in fact isomorphisms.

A ``trivial'' example of a DGA is given by
$(\mathcal{E}^i=\alg{e_1^i,e_2^i},\d)$, where $\deg e_1^i -1 = \deg e_2^i = i$
and $\d e_1^i=e_2^i, \d e_2^i=0$. The main equivalence relation
between DGAs, which we now present, stipulates that adding this
trivial DGA does not change equivalence class.

\begin{definition}
The degree-$i$ \textit{algebraic stabilization} of the DGA
$(\A=\Z\langle a_1,$ $\dots,a_n\rangle,\d)$ is the DGA
$(S^i(\A)=\alg{a_1,\dots,a_n,e_1^i,e_2^i},\d)$, with grading and
differential induced from $\A$ and $\mathcal{E}^i$. Two DGAs
$(\A,\d)$ and $(\tilde{\A},\tilde{\d})$ are \textit{stable tame
isomorphic} if there exist algebraic stabilizations
$(S^{i_1}(\cdots(S^{i_k}(\A))\cdots),\d)$ and
$(S^{\tilde{i}_1}(\cdots(S^{\tilde{i}_{\tilde{k}}}(\tilde{\A}))\cdots),\tilde{\d})$
which are tamely isomorphic.
\end{definition}

Stable tame isomorphism is designed to be a special case of quasi-isomorphism.

\begin{proposition}
If $(\A,\d)$ and $(\tilde{\A},\tilde{\d})$ are stable tame isomorphic,
\label{homology}
then the homologies $H_*(\A,\d)$ and $H_*(\tilde{\A},\tilde{\d})$ are
isomorphic.
\end{proposition}

\noindent The proof of Proposition~\ref{homology} can be found, eg, in
\cite{ENS}.

%*********************************************************************
\subsection{The invariants}
\label{ssec:invariants}

We now introduce the knot and braid invariants. The definition of
the knot DGA depends on the braid definition, which we give first.

Let $B_n$ denote the braid group on $n$ strands, which is generated
by $\sigma_1,\dots,\sigma_{n-1}$, with relations $\sigma_i
\sigma_{i+1} \sigma_i = \sigma_{i+1} \sigma_i \sigma_{i+1}$ for $1
\leq i \leq n-1$, and $\sigma_i \sigma_j = \sigma_j \sigma_i$ for
$|i-j| \geq 2$.  Write $\A_n$ for the tensor algebra on $n(n-1)$
generators which we label as $a_{ij}$ for $1 \leq i,j \leq n$ and
$i\neq j$, and let $\Aut(\A_n)$ denote the group of algebra
automorphisms of $\A_n$.

For each generator $\sigma_k$ of $B_n$, define the tame automorphism
$\hom{\sigma_k} \in \Aut(\A_n)$ by its action on the generators of
$\A_n$:
\[
\hom{\sigma_k}\co\left\{
\begin{array}{ccll}
a_{ki} & \mapsto & -a_{k+1,i} - a_{k+1,k}a_{ki}\phantom{999} & i\neq k,k+1 \\
a_{ik} & \mapsto & -a_{i,k+1} - a_{ik}a_{k,k+1} & i\neq k,k+1 \\
a_{k+1,i} & \mapsto & a_{ki} & i\neq k,k+1 \\
a_{i,k+1} & \mapsto & a_{ik} & i \neq k,k+1 \\
a_{k,k+1} & \mapsto & a_{k+1,k} & \\
a_{k+1,k} & \mapsto & a_{k,k+1} & \\
a_{ij} & \mapsto & a_{ij} & i,j \neq k,k+1.
\end{array}
\right.
\]

\begin{proposition}
$\phi$ extends to a homomorphism from $B_n$ to $\Aut(\A_n)$.
\end{proposition}

\begin{proof}
The fact that $\phi$ preserves the braid relations in $B_n$ is a
direct computation.
\end{proof}

\noindent
We will denote the image of $B\in B_n$ under $\phi$ as $\phi_B \in
\Aut(\A_n)$.

This representation was essentially first studied by Magnus
\cite{Mag} in the context of automorphisms of free groups; more
precisely, he presented a representation on a polynomial algebra in
$\binom{n}{2}$ variables, which we can derive from $\phi$ by
abelianizing and setting $a_{ij} = a_{ji}$ for all $i,j$. Humphries
\cite{Hum} subsequently extended the representation to a polynomial
algebra in $n(n-1)$ variables, and interpreted it in terms of
transvections. Our homomorphism is simply the lift of Humphries'
representation to the corresponding tensor algebra. We note that the
starting point in \cite{II} is a new topological interpretation of
$\phi$ in terms of skein relations.

We can obtain a linear representation from $\phi$ by finding a
finite-dimensional piece of $\A_n$ on which $\phi$ descends; in the
language of \cite{EFM}, this involves linearizing $\phi$ with
respect to some augmentation. Concretely, if we write $\M_n$ as the
subalgebra of $\A_n$ generated by $\{a_{ij}+2\}$, then $\phi$
descends to $\M_n/{\M}_n^2$. If we further quotient by setting
$a_{ij}=a_{ji}$ for all $i,j$, then we obtain an
$\binom{n}{2}$-dimensional representation of $B_n$, which is
precisely the inverse of the Lawrence--Krammer representation
\cite{Kra} with $q=-1$ and $t=1$. Indeed, there is a striking
similarity in appearance between $\phi$ and the Lawrence--Krammer
representation which we further explore in \cite{II}.

Unlike Lawrence--Krammer, $\phi$ does not give a faithful
representation of $B_n$; this is not overly surprising, since $B_n$
has nontrivial center, while $\Aut(\A_n)$ does not. Indeed, a result
from \cite{Hum} states that the kernel of $\phi$ is precisely the
center of $B_n$, which is generated by
$(\sigma_1\cdots\sigma_{n-1})^n$. However, the extension map
$\phi^{\ext}\co B_n \hookrightarrow B_{n+1}
\stackrel{\phi}{\rightarrow} \Aut(\A_{n+1})$ is faithful, a fact
first established algebraically in \cite{Mag}. We will return to
$\phi^{\ext}$ in a moment; geometric proofs of the faithfulness
results can be found in \cite{II}.

We now define the DGA invariant for braids.

\begin{definition}
Let $B\in B_n$ be a braid. Let $\A$ be the graded tensor algebra
on $2n(n-1)$ generators, $\{a_{ij}\,|\,1\leq i,j\leq n,\,i\neq
j\}$ of degree $0$ and $\{b_{ij}\,|\,1\leq i,j\leq n,\,i\neq j\}$
of degree $1$. Define the differential $\d$ on $\A$ by
\[
\d b_{ij} = a_{ij} - \hom{B}(a_{ij}),~~\d a_{ij} =0.
\]
Then we call $(\A,\d)$ the \textit{braid DGA} of $B$.
\label{braiddga}
\end{definition}

\begin{theorem}
If $B,\tilde{B}$ are conjugate in $B_n$, then their braid DGAs are
stable tame isomorphic.
\label{braidthm}
\end{theorem}

\noindent We defer the proof of Theorem~\ref{braidthm} until
Section~\ref{ssec:braidproof}. The DGA invariant for knots is also
derived from the map $\phi$; in fact, we will see in
Section~\ref{ssec:ancillary} that it contains the braid DGA. Before
introducing it, we need some more notation.

As mentioned before, we can define a map $\phi^{\ext}\co
B_n\to\Aut(\A_{n+1})$ induced from $\phi$ through the inclusion $B_n
\hookrightarrow B_{n+1}$ obtained by adding an $(n+1)$-st strand
which does not interact with the other $n$ strands; that is, each
generator $\sigma_i\in B_n$ is mapped to $\sigma_i \in B_{n+1}$ for
$1\leq i\leq n-1$. When dealing with $\phi^{\ext}$, we replace all
indices $(n+1)$ by an asterisk; ie, write $a_{\dott i},a_{i\dott}$
for $a_{n+1,i},a_{i,n+1}$, respectively. Since the $(n+1)$-st strand
does not interact with the other strands, it follows from the
definition of $\phi$ that for $B\in B_n$, $\hom{B}^{\ext}(a_{\dott
i})$ is a linear combination of $a_{\dott j}, 1 \leq j \leq n$, with
coefficients in $\A_n$, and similarly for $\hom{B}^{\ext}(a_{i
\dott})$.

\begin{definition}
The $n\times n$ matrices $\Phil{B}(A)$ and $\Phir{B}(A)$ are defined
by
\[
\hom{B}^{\ext}(a_{i\dott}) = \sum_{j=1}^n (\Phil{B}(A))_{ij} a_{j\dott}
\hspace{0.25in} \textrm{and} \hspace{0.25in}
 \hom{B}^{\ext}(a_{\dott j}) =
\sum_{i=1}^n a_{\dott i} (\Phir{B}(A))_{ij}.
\]
\end{definition}

The matrices $\Phil{B}(A)$ and $\Phir{B}(A)$ can be obtained from
each other by the ``conjugation'' operation discussed in
Section~\ref{ssec:conjugation}. Though $\Phil{B}(A),\Phir{B}(A)$ are
derived from $\phi$, they also determine $\phi$; see
Proposition~\ref{philphir}.

To define the DGA for knots, it is convenient to assemble generators
of the DGA into matrices. For notational purposes, let $A$ denote
the matrix $(a_{ij})$, where we set $a_{ii}=-2$ for all $i$, and
similarly write $B=(b_{ij})$, $C=(c_{ij})$, $D=(d_{ij})$, where
$b_{ij},c_{ij},d_{ij}$ are variables used below with no stipulation
on $b_{ii},c_{ii},d_{ii}$; also, write $\d A$ for the matrix $(\d
a_{ij})$, $\hom{B}(A)$ for $(\hom{B}(a_{ij}))$, and so forth. (There
should hopefully be no confusion between a braid $B$ and the matrix
$B$.) We will sometimes view $A$ as a matrix of variables, so that,
eg, we define $\Phil{B}(M)$ for an $n\times n$ matrix $M$ to be
the evaluation of $\Phil{B}(A)$ when we set $a_{ij}=M_{ij}$ for all
$i,j$. (Note that this requires $M_{ii}=-2$ for all $i$.)

\begin{definition}
Let $B\in B_n$. Write $\A$ for the tensor algebra with the
following generators: $\{a_{ij}\,|\,1\leq i,j\leq n, i\neq j\}$ of
degree $0$; $\{b_{ij},c_{ij}\,|\,1\leq i,j\leq n\}$ of degree $1$;
and $\{d_{ij}\,|\,1\leq i,j\leq n\}$ and $\{e_i\,|\,1\leq i\leq
n\}$ of degree $2$.  Define the differential $\d$ on $\A$ by
\begin{align*}
\d A &= 0 \\
\d B &= (1 - \Phil{B}(A)) \cdot A \\
\d C &= A \cdot (1 - \Phir{B}(A)) \\
\d D &= B \cdot (1 - \Phir{B}(A)) - (1 - \Phil{B}(A)) \cdot C \\
\d e_i &= (B + \Phil{B}(A) \cdot C)_{ii},
\end{align*}
where $\cdot$ denotes matrix multiplication.
\label{knotdga}
Then we call $(\A,\d)$ the \textit{knot DGA} of $B$.
\end{definition}

The fact that $\d^2 = 0$ in the knot DGA is not obvious; while it is
clear that $\d^2 A = \d^2 B = \d^2 C = \d^2 D = 0$, we have $\d^2
e_i = -2 - \left( \Phil{B}(A) \cdot A \cdot \Phir{B}(A)
\right)_{ii}$, and we need Proposition~\ref{matrix} to conclude that
$\d^2 e_i = 0$. The key result for knot DGAs is the following
theorem, which allows the knot DGA to descend from braids to knots.

\begin{theorem}
If the closures of two braids $B,\tilde{B}$ are the same knot,
\label{knotthm}
then the knot DGAs for $B,\tilde{B}$ are stable tame isomorphic.
\end{theorem}

The proof of Theorem~\ref{knotthm} is postponed until
Section~\ref{ssec:fullknotproof}, after we first introduce $HC_0$
and present the easier proof of invariance for $HC_0$.

By Alexander's Theorem, every knot in $\R^3$ can be represented as
a closed braid. Theorem~\ref{knotthm} implies that the following
notion is well-defined.

\begin{definition}
Let $K\subset\R^3$ be a knot. The \textit{knot DGA class} of $K$ is
the equivalence class (under stable tame isomorphism) of the knot
DGA of any braid whose closure is $K$.
\end{definition}

\noindent We will occasionally abuse notation and speak of the
``knot DGA'' of $K$; this is understood to be modulo equivalence.
The knot DGA can equally well be defined for links, for which it
also gives an invariant.

We conclude this section by explicitly computing, as an example, the
knot DGA for the trefoil, which is the closure of $\sigma_1^3\in
B_2$. By the definition of $\phi$, we see that
$\phi^{\ext}_{\sigma_1}$ acts as follows: $a_{12} \mapsto a_{21}$,
$a_{21} \mapsto a_{12}$, $a_{1\dott} \mapsto
-a_{2\dott}-a_{21}a_{1\dott}$, $a_{2\dott} \mapsto a_{1\dott}$,
$a_{\dott 1} \mapsto -a_{\dott 2}-a_{\dott 1}a_{12}$, $a_{\dott 2}
\mapsto a_{\dott 1}$. Iterating $\phi^{\ext}_{\sigma_1}$ three times
on $a_{1\dott},a_{2\dott},a_{\dott 1},a_{\dott 2}$ yields
\[
\Phil{\sigma_1^3}(A) = \left(
\begin{matrix} 2a_{21}-a_{21}a_{12}a_{21} & 1-a_{21}a_{12} \\
-1+a_{12}a_{21} & a_{12}
\end{matrix} \right)\]
$$\Phir{\sigma_1^3}(A) = \left(
\begin{matrix} 2a_{12}-a_{12}a_{21}a_{12} & -1+a_{12}a_{21} \\
1-a_{21}a_{12} & a_{21}
\end{matrix} \right) .\leqno{\rm and}$$
One can now calculate the differential on the knot DGA of
$\sigma_1^3$:
\begin{align*}
\d b_{11} &= -2+3a_{21}-a_{21}a_{12}a_{21}\\
\d c_{11} &= -2+3a_{12}-a_{12}a_{21}a_{12} \\
\d b_{12} &= 2+a_{12}-4a_{21}a_{12}+a_{21}a_{12}a_{21}a_{12}\\
\d c_{21} &= 2+a_{21}-4a_{21}a_{12}+a_{21}a_{12}a_{21}a_{12} \\
\d b_{21} &= -2+a_{21}+a_{12}a_{21}\\
\d c_{12} &= -2+a_{12}+a_{12}a_{21} \\
\d b_{22} &= -2+3a_{12}-a_{12}a_{21}a_{12}\\
\d c_{22} &=-2+3a_{21}-a_{21}a_{12}a_{21}\\
\d d_{11} &= b_{11}-b_{12}-c_{11}+c_{21}+2a_{21}c_{11}-2b_{11}a_{12}-a_{21}a_{12}c_{21}
+b_{12}a_{21}a_{12}\\
& \qquad\qquad -a_{21}a_{12}a_{21}c_{11}+b_{11}a_{12}a_{21}a_{12} \\
\d d_{12} &= b_{11}+b_{12}-c_{12}+c_{22}+2a_{21}c_{12}-b_{12}a_{21}
-a_{21}a_{12}c_{22}\\
& \qquad\qquad -b_{11}a_{12}a_{21}-a_{21}a_{12}a_{21}c_{12} \\
\d d_{21} &= b_{21}-b_{22}-c_{11}-c_{21}+a_{12}c_{21}-2b_{21}a_{12}
+a_{12}a_{21}c_{11}\\
& \qquad\qquad +b_{22}a_{21}a_{12}+b_{21}a_{12}a_{21}a_{12} \\
\d d_{22} &= b_{21}+b_{22}-c_{12}-c_{22}+a_{12}c_{22}-b_{22}a_{21}
+a_{12}a_{21}c_{12}-b_{21}a_{12}a_{21} \\
\d e_1 &= b_{11}+c_{21}+2a_{21}c_{11}-a_{21}a_{12}c_{21}-a_{21}a_{12}a_{21}c_{11} \\
\d e_2 &= b_{22}-c_{12}+a_{12}c_{22}+a_{12}a_{21}c_{12}.
\end{align*}
We will return to this example in future sections.

%*********************************************************************
%*********************************************************************
\section{Motivation from contact geometry}
\label{sec:contact}

Having defined the knot and braid DGAs in the previous section, we
now describe the background in contact geometry leading to the
development of these invariants. It should be noted that work is in
progress to establish rigorously that the DGAs give the relative
contact homology theory which we describe below.

Given a smooth manifold $M$, the cotangent bundle $T^*\!M$ has a
canonical symplectic structure given by $\omega = d\lambda$, where
$\lambda$ is the tautological one-form on $T^*\!M$ sending a tangent
vector in $T^*\!M$ to the pairing between the base point in $T^*\!M$
and the projected tangent vector in $TM$ under the projection
$\pi:T^*\!M\to M$. For any smooth submanifold $N\subset M$, the
conormal bundle $\Lag N = \{\theta\in T^*\!M\,|\, \langle \theta,v
\rangle = 0~\textrm{for all}~v\in T_{\pi(v)}N\}$ over $N$ is then a
Lagrangian submanifold of $T^*\!M$, and smooth isotopy of $N$ leads
to Lagrangian isotopy of $\Lag N$. This setup (with $N$ a knot) was
introduced in the physics literature in \cite{OV}, and was
communicated to the author in this generality by Eliashberg.

Partly because $\Lag N$ is noncompact, it is more convenient for us
to consider a slightly modified setup. Suppose that $M$ is equipped
with a metric. Then the cosphere bundle $ST^*\!M$ of unit covectors
in $T^*\!M$ (``co-oriented contact elements in $M$'') has a
canonical contact structure induced by the one-form $\lambda$. If
$N\subset M$ is a compact smooth submanifold, then the unit conormal
bundle $LN = \Lag N \cap ST^*\!M$ is a compact Legendrian
submanifold of $ST^*\!M$, and smooth isotopy of $N$ leads to
Legendrian isotopy of $LN$.

In recent years, beginning with work of Eliashberg and Hofer,
invariants of Legendrian isotopy in contact manifolds have been
developed, via holomorphic-curve techniques, using relative contact
homology \cite{Eli} and the more general Symplectic Field Theory
\cite{EGH}. In particular, for Legendrian knots in $\R^3$ with the
standard contact structure, Chekanov \cite{Che} developed a purely
combinatorial formulation of relative contact homology; this work
was subsequently extended in \cite{ENS,Ng}. Relative contact
homology was also studied for knots in circle bundles in \cite{Sab},
and for $\R^{2n+1}$ with the standard contact structure in
\cite{EES}. In all of these special cases, relative contact homology
is given by the homology of a certain DGA whose stable tame
isomorphism class is an invariant. With the exception of \cite{EES},
however, the DGAs are defined combinatorially rather than
geometrically, with invariance proofs given by combinatorics as
well. The invariants in the present paper are likewise combinatorial
in nature.

We now specialize, as in \cite{OV}, to the case when $M=\R^3$ and
$K\subset\R^3$ is a knot. Then $LK$ is an embedded Legendrian
$2$-torus in the $5$-dimensional contact manifold $ST^*\R^3 =
\R^3\times S^2$, which is contactomorphic to the $1$-jet space
$J^1(S^2)$ with the standard contact structure. We note that $LK$ is
always topologically unknotted since it has codimension $3$ in
$ST^*\R^3$; however, its Legendrian isotopy type yields information
about the knot $K$. In fact, the knot DGA defined in
Section~\ref{ssec:invariants}, which is a nontrivial invariant,
conjecturally gives the relative contact homology of $LK$ in
$ST^*\R^3$.

A local model suggested by Eliashberg allows us to study conjugacy
classes of braids as a first step. Let $K\subset\R^3$ be any knot
(or, in particular, the unknot), and let $B$ be a braid. A
sufficiently small neighborhood of $LK$ is contactomorphic to
$J^1(LK) = J^1(T^2)$. If we glue $B$ along a tubular neighborhood of
$K$ to produce a new knot $\tilde{K}$, then $L\tilde{K}$ is in an
arbitrarily small neighborhood of $LK$, and thus gives a Legendrian
$2$-torus in $J^1(T^2)$. A smooth isotopy of $B$, or a conjugation
operation on $B$, yields a Legendrian isotopy of this torus in
$J^1(T^2)$. Then the braid DGA from Section~\ref{ssec:invariants}
conjecturally gives the relative contact homology of $L\tilde{K}$ in
$J^1(T^2)$.

Recall from, eg, \cite{Eli} that the relative contact homology for
a Legendrian submanifold $L$ of a contact manifold $(V,\alpha)$ is
the homology of a DGA whose generators are \textit{Reeb chords} of
$L$, ie, flow lines of the Reeb vector field which begin and end
on $L$. The differential is then defined by counting certain
holomorphic maps of punctured disks into the symplectization
$V\times\R$ of $V$, with boundary on $L\times\R$ and punctures
limiting to Reeb chords.

For the case of a knot $K$ in $\R^3$, the contact manifold is $\R^3
\times S^2$, with Reeb vector field pointing in the $\R^3$ fibers,
in the direction specified by the underlying point in $S^2$. It
follows that Reeb chords for $LK$ correspond to ``binormal chords''
to $K$ in $\R^3$, ie, oriented line segments in $\R^3$ with
endpoints on $K$ which are normal to $K$ at both endpoints. For
instance, the unknot given by an ellipse in the plane yields four
Reeb chords, corresponding to the major and minor axes traversed
either way.

If $K$ is the closure of a braid $B\in B_n$, we can embed $B$ in a
neighborhood of an elliptical unknot in such a way that the Reeb
chords for $LK$ come in two families:
\begin{itemize}
\item $2n(n-1)$ ``small chords'' within the neighborhood of the
unknot, four for each pair of braid strands; these divide further
into two Reeb chords where the pair of strands is closest ($a_{ij}$)
and two where the pair is farthest ($b_{ij}$)
\item $4n^2$ ``big chords'' corresponding to the Reeb chords of the
ellipse; each binormal chord of the ellipse yields $n^2$ binormal
chords for $K$ ($c_{ij},d_{ij}$ for the minor-axis chords,
$e_{ij},f_{ij}$ for the major-axis chords), since there are $n$
choices for each endpoint.
\end{itemize}
For the braid itself, the Reeb chords for $LB \subset J^1(T^2)$
correspond to the $2n(n-1)$ ``small chords'' $a_{ij},b_{ij}$.

This gives a rough explanation for the generators of the knot and
braid DGAs; for the knot case, this actually gives the modified knot
DGA defined below in Proposition~\ref{newdga}. Showing that the
differentials of these DGAs actually count the appropriate
holomorphic disks is beyond the scope of this paper. The analytical
details, which use an approach due to Fukaya and Oh \cite{FO} of
counting gradient flow trees, are the subject of work in progress.
As in the theory of Legendrian knots in standard contact $\R^3$,
however, the combinatorial proof of the invariance of our DGAs under
isotopy gives evidence for the validity of our ``computation'' of
relative contact homology.

%*********************************************************************
%*********************************************************************
\section{Knot and braid contact homology}
\label{sec:HC}

In this section, we examine the homology of the DGAs defined in
Section~\ref{sec:defs}, and focus on a particular piece for which
isotopy invariance is easy to establish. We will return to
invariance proofs of the full DGAs in Section~\ref{sec:fullproof}.

%*********************************************************************
\subsection{Definitions}
\label{ssec:HCdef}

In order to use the knot and braid DGAs as invariants, we need to
know when two DGAs are stable tame isomorphic. This can be done
through ``computable'' invariants of equivalence classes of DGAs,
including Poincar\'e polynomials \cite{Che} and the characteristic
algebra \cite{Ng}.

In the case of braid and knot DGAs, it turns out that the first
order Poincar\'e polynomial over $\Z_2$ corresponding to the trivial
augmentation gives no interesting invariant of the braid conjugacy
class or the knot. (See however Section~\ref{ssec:T2} for a
different result for a nontrivial augmentation over $\Z_2$, and
Section~\ref{ssec:lin} for linearization over $\Z$.) This is because
to first order over $\Z_2$, the homomorphism $\phi$ factors through
the projection $B_n \to S_n$ to the symmetric group. However, the
degree $0$ part of the characteristic algebra already yields
interesting invariants; we recall from \cite{Ng} that the
characteristic algebra of a DGA $(\A,\d)$ is the quotient of $\A$ by
the two-sided ideal generated by $\im\d$. Before proceeding further,
we define the contact homology of a braid or knot.

\begin{definition}
The \textit{contact homology} of a braid $B$, written $HC_*(B)$, is
defined to be the (graded) homology of the braid DGA of $B$. The
\textit{contact homology} of a knot $K$, written $HC_*(K)$, is the
(graded) homology of the knot DGA of any braid whose closure is $K$.
More generally, if $R$ is any ring, then we define as usual
$HC_*(B;R)$ and $HC_*(K;R)$ to be the homology of the appropriate
DGA tensored with $R$.
\end{definition}

\begin{proposition}
The degree $0$ contact homology of knots and braids is the degree
$0$ part of the characteristic algebra of the corresponding DGA.
More explicitly, for $B\in B_n$, let $\I^{\braid}_B$ be the
two-sided ideal of $\A_n$ generated by $\{a_{ij}-\hom{B}(a_{ij})\}$,
while $\I^{\knot}_B$ is the two-sided ideal of $\A_n$ generated by
the entries of the two matrices $(1 - \Phil{B}(A))\cdot A$ and
$A\cdot (1 - \Phir{B}(A))$. Then $HC_0(B) = \A_n/\I^{\braid}_B$ and
$HC_0(K) = \A_n/\I^{\knot}_B$, where $K$ is the closure of $B$.
\label{charalg}
\end{proposition}

\begin{proof}
Since the braid and knot DGAs of $B\in B_n$ have no generators in
negative degree, the entirety of $\A_n$ (which is the degree $0$
part of the algebra) consists of cycles. On the other hand, the
image of $\d$ in degree $0$ is precisely the ideal generated by the
images under $\d$ of the generators in degree $1$, and so the ideal
of boundaries in degree $0$ is given by $\I^{\braid}_B$,
$\I^{\knot}_B$ for the braid and knot DGA, respectively.
\end{proof}

For example, for the trefoil, $HC_0(3_1)$ is given by $\Z\langle
a_{12},a_{21}\rangle$ modulo the relations provided by $\d b_{ij},\d
c_{ij}$ as calculated in Section~\ref{ssec:invariants}. From $\d
b_{21}$ and $\d c_{12}$, we have relations $-2+a_{21}+a_{12}a_{21} =
-2+a_{12}+a_{12}a_{21} = 0$, and hence $a_{12}=a_{21}$. If we set
$a_{12}=a_{21}=x$, then the $\d b_{ij},\d c_{ij}$ relations become
$\{x^3-3x+2, x^4-4x^2+x+2, x^2+x-2\}$. The gcd of these polynomials
is $x^2+x-2$, and so $HC_0(3_1) \cong \Z[x]/(x^2+x-2)$.

In practice, we will study $HC_0$ rather than the full braid and
knot DGAs, because it is relatively easy to compute, as the trefoil
example demonstrates. It is immediate from
Proposition~\ref{homology} and Theorems~\ref{braiddga} and
\ref{knotdga} that $HC_0$ gives an invariant of braid group
conjugacy classes and knots. However, there is a direct proof of the
invariance of $HC_0$ which is much simpler notationally than the
proofs of Theorems~\ref{braiddga} and \ref{knotdga}, while
containing the main ideas from these proofs. We will give this
direct proof below for braids, and in
Section~\ref{ssec:HCknotinvarianceproof} for knots.

\begin{theorem}
Up to isomorphism,
\label{HCbraidinvariance}
$HC_0(B)$ is an invariant of the conjugacy class of the braid $B$.
\end{theorem}

\begin{proof}
Suppose that $B,\tilde{B}$ are conjugate in the braid group $B_n$,
so that $\tilde{B} = C^{-1}BC$ for some $C\in B_n$. We use the
notation of Proposition~\ref{charalg}.

First note that $\I^{\braid}_B \subset \A_n$ is the ideal generated
by the image of the map $1-\hom{B}\co \A_n\to\A_n$: if
$x_1,x_2\in\A_n$, then $(1-\hom{B})(x_1x_2) = (1-\hom{B})(x_1) x_2 +
\hom{B}(x_1) (1-\hom{B})(x_2)$ is in the ideal generated by
$(1-\hom{B})(x_1)$ and $(1-\hom{B})(x_2)$. Similarly,
$\I^{\textrm{braid}}_{\tilde{B}}$ is the ideal generated by the
image of $1-\hom{\tilde{B}} = 1 - \hom{C}^{-1} \hom{B} \hom{C}$. It
follows by inspection of the commutative diagram
\[
\xymatrix{
\A_n \ar[rr]^{1-\hom{B}} \ar[d]_{\hom{C}}^{\cong} && \A_n \ar[d]^{\hom{C}}_{\cong} \\
\A_n \ar[rr]^{1-\hom{C}^{-1} \hom{B} \hom{C}} && \A_n
}
\]
that $\I^{\braid}_{\tilde{B}}$ is the image of $\I^{\braid}_B$ under
the automorphism $\hom{C}$ on $\A_n$, and hence that $HC_0(B) =
\A_n/\I^{\braid}_B$ and $HC_0(\tilde{B}) =
\A_n/\I^{\braid}_{\tilde{B}}$ are isomorphic.
\end{proof}

The corresponding invariance result for knots requires a few
preparatory results, which are interesting in their own right; we
give these in the following section, and prove $HC_0$ invariance for
knots in Section~\ref{ssec:HCknotinvarianceproof}.

%*********************************************************************
\subsection{Ancillary results}
\label{ssec:ancillary}

Here we collect a number of results about $\phi$, $\Phil{}$, and $\Phir{}$
which will be crucial for future sections.

\begin{proposition}[Chain Rule]
If $B_1,B_2\in B_n$, then
\label{chainrule}
\[
\Phi^L_{B_1 B_2} (A) = \Phil{B_2}(\hom{B_1}(A)) \cdot \Phil{B_1}(A)
\quad \textrm{and} \quad \Phir{B_1 B_2}(A) = \Phir{B_1}(A)
\cdot \Phir{B_2} (\hom{B_1}(A)).
\]
\end{proposition}

\begin{proof}
We will establish the result for $\Phi^L$; the proof for $\Phi^R$
is completely analogous. By the definition of $\Phi^L$, we have
\[
\hom{B_1B_2}^{\ext}(a_{i \dott}) = \sum_{j=1}^n \left(
\Phil{B_1B_2}(A) \right)_{ij} a_{j \dott};
\]
on the other hand,
\begin{align*}
\hom{B_1B_2}^{\ext}(a_{i\dott}) &= \hom{B_1}^{\ext} \hom{B_2}^{\ext}
(a_{i\dott})
\\
&= \sum_{j=1}^n \hom{B_1}^{\ext} \left( \left( \Phil{B_2}(A)
\right)_{ij} a_{j\dott}\right)
\\
&= \sum_{j=1}^n \left(\Phil{B_2}(\hom{B_1}(A))\right)_{ij}
\hom{B_1}^{\ext}(a_{j\dott})
\\
&= \sum_{j,\ell=1}^n \left(\Phil{B_2}(\hom{B_1}(A))\right)_{ij}
(\Phil{B_1}(A))_{j\ell} a_{\ell\dott},
\end{align*}
and the result follows.

A word of explanation is in order for the penultimate equality in
the chain above, because it may seem backwards. Applying
$\hom{B_1}^{\ext}$ to $\left(\Phil{B_2}(A)\right)_{ij}a_{j\dott}$
entails replacing each $a_{ij}$ by $\hom{B_1}(a_{ij})$ and each
$a_{j\dott}$ by $\hom{B_1}^{\ext}(a_{j\dott})$, and the result is
$\left(\Phil{B_2}(\hom{B_1}(A))\right)_{ij} \hom{B_1}^{\ext}(a_{j\dott})$.
\end{proof}

\begin{corollary}
For any braid $B$, \label{inverse} $$\Phil{B^{-1}}(\hom{B}(A)) =
(\Phil{B}(A))^{-1}\quad\text{and}\quad\Phir{B^{-1}}(\hom{B}(A)) =
(\Phir{B}(A))^{-1}.$$
\end{corollary}

\begin{proof}
Set $B_1=B$, $B_2=B^{-1}$ in Proposition~\ref{chainrule}, and use
the fact that $\Phil{\textrm{id}}(A) = \Phir{\textrm{id}}(A) = 1.$
\end{proof}

\begin{lemma}
The matrices $\Phil{\sigma_k}(A)$ and $\Phir{\sigma_k}(A)$ are
identical to the $n\times n$ identity matrix, outside of the
$2\times 2$ submatrix formed by rows $k,k+1$ and columns $k,k+1$,
which is given by
$\left( \begin{smallmatrix} -a_{k+1,k} & -1 \\ 1 & 0 \end{smallmatrix} \right)$
for $\Phil{\sigma_k}(A)$ and
$\left( \begin{smallmatrix} -a_{k,k+1} & 1 \\ -1 & 0 \end{smallmatrix} \right)$
for $\Phir{\sigma_k}(A)$.
\label{sigmak}
\end{lemma}

\begin{proof}
This is immediate from the definitions of $\phi$, $\Phil{}$, and $\Phir{}$.
\end{proof}

\begin{proposition}
For any braid $B$,
$\hom{B}(A) = \Phil{B}(A) \cdot A \cdot \Phir{B}(A).$
\label{matrix}
\end{proposition}

\begin{proof}
When $B=\sigma_k$ for some $k$, the identity can be verified by direct
calculation, using Lemma~\ref{sigmak}.
We conclude that the identity also holds for $B=\sigma_k^{-1}$,
either by another direct calculation, or by using the result for
$B=\sigma_k$: since
$\hom{\sigma_k}(A) = \Phil{\sigma_k}(A) \cdot A \cdot \Phir{\sigma_k}(A),$
we have
\[
A = \Phil{\sigma_k}(\hom{\sigma_k}^{-1}(A)) \cdot \hom{\sigma_k}^{-1}(A) \cdot
\Phir{\sigma_k}(\hom{\sigma_k}^{-1}(A)) =
\left( \Phil{\sigma_k^{-1}}(A) \right)^{-1}\!\! \cdot \hom{\sigma_k^{-1}}(A) \cdot
\left( \Phir{\sigma_k^{-1}}(A) \right)^{-1},
\]
where the last equality comes from Corollary~\ref{inverse};
the desired identity for $B=\sigma_k^{-1}$
follows.

We now show that if the identity holds for $B=B_1$ and $B=B_2$, then
it holds for $B=B_1B_2$; the proposition then follows by induction.
Indeed, assuming the identity for $B=B_1$ and $B=B_2$, we have
\begin{align*}
\hom{B_1B_2}(A) &= \hom{B_1}(\hom{B_2}(A)) \\
&=
\Phil{B_2}(\hom{B_1}(A)) \cdot \hom{B_1}(A) \cdot \Phir{B_2}(\hom{B_1}(A)) \\
&= \Phil{B_2}(\hom{B_1}(A)) \cdot \Phil{B_1}(A) \cdot A \cdot
\Phir{B_1}(A)
\cdot \Phir{B_2}(\hom{B_1}(A)) \\
&= \Phil{B_1B_2}(A) \cdot A \cdot \Phir{B_1B_2}(A),
\end{align*}
where we use the chain rule for the final equality.
\end{proof}

A more geometric proof of Proposition~\ref{matrix} is given in
\cite{II}. As remarked in Section~\ref{ssec:invariants},
Proposition~\ref{matrix} implies that the knot DGA does in fact
satisfy $\d^2=0$. It also follows from Proposition~\ref{matrix} that
the full homomorphism $\hom{B}$ can be deduced from the matrices
$\Phil{B}$ and $\Phir{B}$, which in turn can be deduced from each
other via the conjugation map to be discussed in
Section~\ref{ssec:conjugation}. Thus we may view $\Phil{B}$ and
$\Phir{B}$ as ``square roots'' of $\hom{B}$.

One final result we will need in the future is a reformulation of
the knot DGA in terms of an equivalent DGA which has more generators
but a slightly simpler definition.

\begin{proposition}
Let $B\in B_n$. The knot DGA for $B$ is stable tame isomorphic to
the ``modified knot DGA'' which has generators $\{a_{ij}\,|\,1\leq
i,j\leq n, i\neq j\}$ of degree $0$, $\{b_{ij}\,|\,1\leq i,j\leq n,
i\neq j\}$ and $\{c_{ij},d_{ij}\,|\,1\leq i,j\leq n\}$ of degree
$1$, and $\{e_{ij},f_{ij}\,|\,1\leq i,j\leq n\}$ of degree $2$, with
differential \label{newdga}
\begin{align*}
\d A &= 0 \\
\d B &= A - \hom{B}(A) \\
\d C &= (1 - \Phil{B}(A)) \cdot A \\
\d D &= A \cdot (1 - \Phir{B}(A)) \\
\d E &= B - D - C \cdot \Phir{B}(A) \\
\d F &= B - C - \Phil{B}(A) \cdot D.
\end{align*}
Here we set $b_{ii}=0$ for all $i$.
\end{proposition}

\begin{proof}
In the modified knot DGA, replace $e_{ij}$ by $e_{ij}+f_{ij}$ for
all $i,j$; after applying the induced tame automorphism, we obtain
the same differential, except with $\d E = C\cdot
(1-\Phir{B}(A))-(1-\Phil{B}(A))\cdot D$. Now replace $b_{ij}$ by
$(B+C+\Phil{B}(A)\cdot D)_{ij}$ for $i\neq j$, which constitutes
another tame automorphism; the differential is again unchanged,
except that $\d B = 0$ (this follows from Proposition~\ref{matrix})
and $\d F = B - \operatorname{diag}(C+\Phil{B}(A)\cdot D)$, where
$\operatorname{diag}$ replaces all non-diagonal entries by $0$. Now
we have $\d f_{ij} = b_{ij}$ for $i\neq j$, and hence, up to
stabilization, we may drop the generators $\{b_{ij},f_{ij}\,|\,i\neq
j\}$ from the DGA. The resulting DGA is precisely the knot DGA for
$B$, after we relabel $c_{ij},d_{ij},e_{ij},f_{ii}$ by
$b_{ij},c_{ij},d_{ij},-e_i$, respectively.
\end{proof}

An important consequence of Proposition~\ref{newdga} is that the
knot DGA for a braid can be obtained from the braid DGA, up to
stable tame isomorphism, by adding generators and appropriately
extending the definition of the differential.

\begin{corollary}
If $K$ is the closure of $B$, then $HC_0(K)$ is a quotient of $HC_0(B)$.
\label{quotientcor}
\end{corollary}

\begin{proof}
Immediate from Propositions~\ref{homology} and~\ref{newdga}.
\end{proof}

%*********************************************************************
\subsection{Invariance of $HC_0$ for knots}
\label{ssec:HCknotinvarianceproof}

In this section, we prove the following result.

\begin{theorem}
Up to isomorphism, $HC_0(K)$ is a well-defined invariant of a knot $K$.
\label{HCknotinvariance}
\end{theorem}

\begin{proof}
If $B$ is a braid, we temporarily write $\HCknot(B)$ for the
$0$-dimensional homology of the knot DGA of $B$, as expressed in
Proposition~\ref{charalg}. We wish to show that if $B,\tilde{B}$ are
two braids whose closure is the same knot, then
$\HCknot(B),\HCknot(\tilde{B})$ are isomorphic. By Markov's Theorem,
it suffices to check this when $B,\tilde{B}$ are related by one of
three operations: conjugation, positive stabilization, and negative
stabilization.

\medskip
\textbf{Conjugation}\qua Here $B,\tilde{B} \in B_n$ satisfy
$\tilde{B} = C^{-1} B C$ for some $C\in B_n$.

Write $\HCknot(B) = \A_n/\I$, where $\I$ is generated by the entries
of $(1-\Phil{B}(A))\cdot A$ and $A\cdot (1-\Phir{B}(A))$. For
clarity, write $\HCknot(\tilde{B}) = \tilde{\A}_n/\tilde{\I}$, with
$\tilde{\A}_n$ generated by $\{\tilde{a}_{ij}\}$ and $\tilde{\I}$
generated by the entries of $(1+\Phil{\tilde{B}}(\tilde{A}))\cdot
\tilde{A}$ and $\tilde{A}\cdot (1+\Phir{\tilde{B}}(\tilde{A}))$,
where $\tilde{A} = (\tilde{a}_{ij})$ and $\tilde{a}_{ii}=-2$.

We claim that the identification $\tilde{A} = \hom{C}(A)$, ie,
$\tilde{a}_{ij} = \hom{C}(a_{ij})$ for all $i,j$, maps $\tilde{\I}$
into $\I$. This would suffice to prove that $\HCknot(B) \cong
\HCknot(\tilde{B})$, since a symmetrical argument with $B$ and
$\tilde{B}$ interchanged shows that the same identification in the
other direction also maps $\I$ into $\tilde{\I}$, and hence $\I$ and
$\tilde{\I}$ coincide.

Assume that $\tilde{A} = \hom{C}(A)$.
From Proposition~\ref{chainrule} and Corollary~\ref{inverse}, we have
\[
\Phil{\tilde{B}}(\tilde{A}) =
\Phil{BC}(A) \cdot \Phil{C^{-1}}(\hom{C}(A))) =
\Phil{C}(\hom{B}(A)) \cdot \Phil{B}(A) \cdot (\Phil{C}(A))^{-1},
\]
and similarly
$\Phir{\tilde{B}}(\tilde{A}) =
(\Phir{C}(A))^{-1} \cdot \Phir{B}(A) \cdot \Phir{C}(\hom{B}(A))$.
It follows by Proposition~\ref{matrix} that
\begin{align*}
(1-\Phil{\tilde{B}}(\tilde{A})) \cdot \tilde{A}
&= \left(1 -
\Phil{C}(\hom{B}(A)) \cdot \Phil{B}(A) \cdot
(\Phil{C}(A))^{-1}\right)
\cdot \Phil{C}(A) \cdot A \cdot \Phir{C}(A) \\
&= \Phil{C}(A) \cdot (1 - \Phil{B}(A))\cdot A \cdot \Phir{C}(A) \\
& \qquad - \left( \Phil{C}(\hom{B}(A)) - \Phil{C}(A) \right) \cdot
\Phil{B}(A) \cdot A \cdot \Phir{C}(A).
\end{align*}
Now the entries of $(1-\Phil{B}(A)) \cdot A$ are in $\I$ by
definition; on the other hand, the entries of
\[
\hom{B}(A) - A = \Phil{B}(A) \cdot A \cdot \Phir{B}(A) - A
= -(1-\Phil{B}(A))\cdot A \cdot \Phir{B}(A) - A \cdot (1-\Phir{B}(A))
\]
are in $\I$, and thus so are the entries of $\Phil{C}(\hom{B}(A)) -
\Phil{C}(A)$. Hence the entries of $(1-\Phil{\tilde{B}}(\tilde{A}))
\cdot \tilde{A}$ are in $\I$. Similarly, the entries of $\tilde{A}
\cdot (1-\Phir{\tilde{B}}(\tilde{A}))$ are also in $\I$, and
therefore $\tilde{\I} \subset \I$, as desired.

\medskip
\textbf{Positive stabilization}\qua Here $\tilde{B}$ is
obtained from $B$ by adding an extra strand, which we label $0$, to
the braid, and setting $\tilde{B} = B\sigma_0$.

Analogously to before, write $\HCknot(B) = \A_n/\I$ and
$\HCknot(\tilde{B}) = \tilde{\A}/\tilde{\I}$, where now $\tilde{\A}$
is obtained from $\A_n$ by adding generators $a_{i0},a_{0i}$ for $1
\leq i \leq n$. Abbreviate $\Phil{B}(A),\Phir{B}(A)$ by
$\Phi^L,\Phi^R$. By using Proposition~\ref{chainrule} and
Lemma~\ref{sigmak}, we can easily compute the matrices for
$\Phil{B\sigma_0},\Phir{B\sigma_0}$ in terms of $\Phi^L,\Phi^R$:
\begin{xalignat*}{2}
\left(\Phil{B\sigma_0}\right)_{00} &= -\hom{B}(a_{10}) & \qquad
\left(\Phir{B\sigma_0}\right)_{00} &= -\hom{B}(a_{01}) \\
\left(\Phil{B\sigma_0}\right)_{0i} &= -\Phi^L_{1i} & \qquad
\left(\Phir{B\sigma_0}\right)_{i0} &= -\Phi^R_{i1} \\
\left(\Phil{B\sigma_0}\right)_{10} &= 1 & \qquad
\left(\Phir{B\sigma_0}\right)_{01} &= 1 \\
\left(\Phil{B\sigma_0}\right)_{1i} &= 0 & \qquad
\left(\Phir{B\sigma_0}\right)_{i1} &= 0 \\
\left(\Phil{B\sigma_0}\right)_{j0} &= 0 & \qquad
\left(\Phir{B\sigma_0}\right)_{0j} &= 0 \\
\left(\Phil{B\sigma_0}\right)_{ji} &= \Phi^L_{ji} & \qquad
\left(\Phir{B\sigma_0}\right)_{ij} &= \Phi^R_{ij}.
\end{xalignat*}
Here $i,j$ are any indices such that $i\geq 1$ and $j\geq 2$. It is
then straightforward to calculate the generators of $\tilde{\I}$;
they are the entries of $(1-\Phil{\tilde{B}}(A)) \cdot A$ and $A
\cdot (1-\Phir{\tilde{B}}(A))$, which we write as
$\tilde{\d}\tilde{B}$ and $\tilde{\d}\tilde{C}$ for future use in
Section~\ref{ssec:fullknotproof}:
\begin{xalignat*}{2}
\tilde{\d} \tilde{b}_{00} &= -2-\hom{B}(a_{10}) & \qquad
\tilde{\d} \tilde{c}_{00} &= -2-\hom{B}(a_{01}) \\
\tilde{\d} \tilde{b}_{0i} &= a_{0i} + \hom{B}(a_{10})a_{0i} + \Phi^L_{1\ell}a_{\ell i} & \qquad
\tilde{\d} \tilde{c}_{i0} &= a_{i0} + a_{i0}\hom{B}(a_{01}) + a_{i\ell}\Phi^R_{\ell 1} \\
\tilde{\d} \tilde{b}_{10} &= a_{10}+2 & \qquad
\tilde{\d} \tilde{c}_{01} &= a_{01}+2  \\
\tilde{\d} \tilde{b}_{1i} &= a_{1i} - a_{0i}  & \qquad
\tilde{\d} \tilde{c}_{i1} &= a_{i1} - a_{i0}  \\
\tilde{\d} \tilde{b}_{j0} &= a_{j0} - \Phi^L_{j\ell} a_{\ell 0}  & \qquad
\tilde{\d} \tilde{c}_{0j} &= a_{0j} - a_{0\ell} \Phi^R_{\ell j}  \\
\tilde{\d} \tilde{b}_{ji} &= a_{ji} - \Phi^L_{j\ell} a_{\ell i}  & \qquad
\tilde{\d} \tilde{c}_{ij} &= a_{ij} - a_{i\ell} \Phi^R_{\ell j},
\end{xalignat*}
where, as before, $i\geq 1$ and $j\geq 2$, and any monomial
involving the index $\ell$ is understood to be summed from $\ell=1$
to $\ell=n$. (To obtain the expressions for $\tilde{b}_{00}$ and
$\tilde{c}_{00}$, we use the facts that
$\hom{B}(a_{10})=\Phi^L_{1\ell}a_{\ell 0}$ and
$\hom{B}(a_{01})=a_{0\ell}\Phi^R_{\ell 1}$, which follow from the
definitions of $\Phi^L$ and $\Phi^R$.)

In $\tilde{\A}/\tilde{\I}$, we thus have $a_{0i} = a_{1i}$ and
$a_{i0} = a_{i1}$ for $i \geq 1$; in particular, $a_{01}=a_{10}=-2$.
Using these relations to replace all generators of the form $a_{0i}$
and $a_{i0}$, we find that the remaining relations in
$\tilde{\A}/\tilde{\I}$ give precisely the generators of $\I$,
namely $a_{ij} - \Phi^L_{i\ell} a_{\ell j}$ and $a_{ij} - a_{i\ell}
\Phi^R_{\ell j}$ for $1\leq i,j\leq n$. Thus $\A_n/\I \cong
\tilde{\A}/\tilde{\I}$, as desired.

\medskip
\textbf{Negative stabilization}\qua Here $\tilde{B}$ is
obtained from $B$ by adding an extra strand, which we label $0$, to
the braid, and setting $\tilde{B} = B\sigma_0^{-1}$.

This case is very similar to the case of positive stabilization,
with a slightly different computation. Use the same notation as for
positive stabilization; then we have
\begin{xalignat*}{2}
\left(\Phil{B\sigma_0}\right)_{00} &= 0 & \qquad
\left(\Phir{B\sigma_0}\right)_{00} &= 0 \\
\left(\Phil{B\sigma_0}\right)_{0i} &= \Phi^L_{1i} & \qquad
\left(\Phir{B\sigma_0}\right)_{i0} &= \Phi^R_{i1} \\
\left(\Phil{B\sigma_0}\right)_{10} &= -1 & \qquad
\left(\Phir{B\sigma_0}\right)_{01} &= -1 \\
\left(\Phil{B\sigma_0}\right)_{1i} &= -\hom{B}(a_{01}) \Phi^L_{1i} & \qquad
\left(\Phir{B\sigma_0}\right)_{i1} &= -\Phi^R_{i1} \hom{B}(a_{10}) \\
\left(\Phil{B\sigma_0}\right)_{j0} &= 0 & \qquad
\left(\Phir{B\sigma_0}\right)_{0j} &= 0 \\
\left(\Phil{B\sigma_0}\right)_{ji} &= \Phi^L_{ji} & \qquad
\left(\Phir{B\sigma_0}\right)_{ij} &= \Phi^R_{ij}
\end{xalignat*}
for $i\geq 1$ and $j\geq 2$, and so the generators of $\tilde{\I}$
are
\begin{xalignat*}{2}
\tilde{\d} \tilde{b}_{00} &= -2-\hom{B}(a_{10}) & \qquad
\tilde{\d} \tilde{c}_{00} &= -2-\hom{B}(a_{01}) \\
\tilde{\d} \tilde{b}_{0i} &= a_{0i} - \Phi^L_{1\ell} a_{\ell i} & \qquad
\tilde{\d} \tilde{c}_{i0} &= a_{i0} - a_{i\ell} \Phi^R_{\ell 1} \\
\tilde{\d} \tilde{b}_{10} &= a_{10} - 2 + \hom{B}(a_{01}) \hom{B}(a_{10}) & \qquad
\tilde{\d} \tilde{c}_{01} &= a_{01} -2 + \hom{B}(a_{01}) \hom{B}(a_{10}) \\
\tilde{\d} \tilde{b}_{1i} &= a_{1i} + a_{0i} + \hom{B}(a_{01}) \Phi^L_{1\ell} a_{\ell i} & \qquad
\tilde{\d} \tilde{c}_{i1} &= a_{i1} + a_{i0} + a_{i\ell} \Phi^R_{\ell 1} \hom{B}(a_{10}) \\
\tilde{\d} \tilde{b}_{j0} &= a_{j0} - \Phi^L_{j\ell} a_{\ell 0}  & \qquad
\tilde{\d} \tilde{c}_{0j} &= a_{0j} - a_{0\ell} \Phi^R_{\ell j}  \\
\tilde{\d} \tilde{b}_{ji} &= a_{ji} - \Phi^L_{j\ell} a_{\ell i}  & \qquad
\tilde{\d} \tilde{c}_{ij} &= a_{ij} - a_{i\ell} \Phi^R_{\ell j}.
\end{xalignat*}
In $\tilde{\A}/\tilde{\I}$, we have $a_{0i} = \Phi^L_{1\ell} a_{\ell
i}$ and $a_{i0} = a_{i\ell} \Phi^R_{\ell 1}$. When we use these
relations to eliminate $a_{0i},a_{i0}$, we find that the remaining
relations in $\tilde{\A}/\tilde{\I}$ are precisely the relations in
$\A_n/\I$, as before. Hence $\tilde{\A}/\tilde{\I} \cong \A_n/\I$.

We have shown that $\HCknot(B)$ is invariant under the Markov moves, and hence
it gives a well-defined knot invariant, as desired.
\end{proof}

%*********************************************************************
%*********************************************************************
\section{Invariance proofs for braid and knot DGAs}
\label{sec:fullproof}

In this section, we prove Theorems~\ref{braidthm} and \ref{knotthm},
the invariance results for braid and knot DGAs. As mentioned before,
the proofs are essentially more involved versions of the
corresponding proofs for $HC_0$ given in Section~\ref{sec:HC}.

%*********************************************************************
\subsection{Proof of Theorem~\ref{braidthm}}
\label{ssec:braidproof}

It suffices to establish Theorem~\ref{braidthm} under the
assumption that $\tilde{B} = \sigma_k^{-1} B \sigma_k$ for some
$k$. Then the braid DGAs for $B$ and $\tilde{B}$ are generated by
$a_{ij},b_{ij}$ and $\tilde{a}_{ij},\tilde{b}_{ij}$, respectively,
for $1\leq i\neq j\leq n$.  If we abbreviate $\hom{B}$ by $\phi$
and $\hom{\sigma_k}$ by $\phi_k$, then the differentials for the
braid DGAs are given by $\d b_{ij} = a_{ij} - \phi(a_{ij})$ and
$\tilde{\d} \tilde{b}_{ij} = \tilde{a}_{ij} - \phi_k^{-1} \phi
\phi_k (\tilde{a}_{ij})$.  (By abusing notation, we will consider
$\phi$ and $\phi_k$ to act either on the algebra generated by the
$a_{ij}$ or on the algebra generated by the $\tilde{a}_{ij}$.) We
will find a tame isomorphism between the algebras which sends
one differential to the other.

Set $\tilde{a}_{ij} = \phi_k (a_{ij})$ for all $i,j$. Then
\[
\tilde{\d} \tilde{b}_{ij} = \tilde{a}_{ij} - \phi_k^{-1} \phi
\phi_k (\tilde{a}_{ij}) = \phi_k (a_{ij}) - \phi \phi_k (a_{ij}).
\]
Now make the following identifications:
\begin{alignat*}{4}
\tilde{b}_{ki} &= - b_{k+1,i} - \phi(a_{k+1,k})b_{ki} - b_{k+1,k} a_{ki} & \qquad & i\neq
k,k+1 \\
\tilde{b}_{ik} &= - b_{i,k+1} - b_{ik} \phi(a_{k,k+1}) - a_{ik} b_{k,k+1} & \qquad & i\neq
k,k+1 \\
\tilde{b}_{k+1,i} &= b_{ki} & \qquad & i\neq k,k+1 \\
\tilde{b}_{i,k+1} &= b_{ik} & \qquad & i\neq k,k+1 \\
\tilde{b}_{k,k+1} &= b_{k+1,k} & \\
\tilde{b}_{k+1,k} &= b_{k,k+1} & \\
\tilde{b}_{ij} &= b_{ij} & \qquad & i,j\neq k,k+1.
\end{alignat*}
It is straightforward to check from $\d b_{ij} = a_{ij} -
\phi(a_{ij})$ that $\d \tilde{b}_{ij} = \phi_k (a_{ij}) - \phi
\phi_k (a_{ij})$ for all $i,j$.  For instance, for $i\neq k,k+1$, we
have
\begin{align*}
\d \tilde{b}_{ik} &= - \d b_{i,k+1} - (\d b_{ik}) \phi(a_{k,k+1}) -
a_{ik} (\d b_{k,k+1}) \\
&= - a_{i,k+1} + \phi(a_{i,k+1}) + \phi(a_{ik})
\phi(a_{k,k+1}) - a_{ik}a_{k,k+1} \\
&= \phi_k (a_{ik}) - \phi \phi_k (a_{ik}).
\end{align*}

Using the above definitions for $\tilde{a}_{ij},\tilde{b}_{ij}$ in
terms of $a_{ij},b_{ij}$, we conclude that $\d = \tilde{\d}$, and so
the map sending $a_{ij}$ to $\tilde{a}_{ij}$ and $b_{ij}$ to
$\tilde{b}_{ij}$ sends $\d$ to $\tilde{\d}$. On the other hand, it
is easy to check that this map is tame; the map on the $a_{ij}$ is
given by the tame isomorphism $\phi_k$, while the map on the
$b_{ij}$ is tame by inspection of the above definition for
$\tilde{b}_{ij}$ in terms of $b_{ij}$. Thus the braid DGAs for $B$
and $\sigma_k^{-1} B \sigma_k$ are tamely isomorphic, as desired.
\qed

%*********************************************************************
\subsection{Proof of Theorem~\ref{knotthm}}
\label{ssec:fullknotproof}

As in the proof of Theorem~\ref{HCknotinvariance}, it suffices to
show that the equivalence class of the knot DGA of a braid is
invariant under the Markov moves.

\medskip
\textbf{Conjugation}\qua Let $B\in B_n$ be a braid, and let
$\tilde{B} = \sigma_k^{-1} B \sigma_k$; we show that the modified
knot DGAs (see Proposition~\ref{newdga}) for $B$ and $\tilde{B}$ are
tamely isomorphic.

For clarity, we distinguish between the modified knot DGAs for $B$
and $\tilde{B}$ by using tildes on the generators and differential
of the knot DGA for $\tilde{B}$. As in the proof of
Theorem~\ref{braidthm}, we will exhibit an identification between
the two sets of generators so that $\d=\tilde{\d}$.

By Lemma~\ref{sigmak}, the matrices
$\Phil{\sigma_k}(A)-\Phil{\sigma_k}(\hom{B}(A))$ and
$\Phir{\sigma_k}(A)-\Phir{\sigma_k}(\hom{B}(A))$ are both
identically zero except in the $(k,k)$ entry, where they are
$-a_{k+1,k}+\hom{B}(a_{k+1,k})$ and $-a_{k,k+1}+\hom{B}(a_{k,k+1})$,
respectively. Hence if we define $n\times n$ matrices
$\Theta^L_k(A)$ and $\Theta^R_k(A)$ to be zero except in the $(k,k)$
entry, where they are $-b_{k+1,k}$ and $-b_{k,k+1}$, respectively,
then we have
\[
\d \Theta^L_k(A) = \Phil{\sigma_k}(A)-\Phil{\sigma_k}(\hom{B}(A))
\quad \textrm{and} \quad
\d \Theta^R_k(A) = \Phir{\sigma_k}(A)-\Phir{\sigma_k}(\hom{B}(A)).
\]

Identify the $\tilde{a}_{ij},\tilde{b}_{ij}$ with the
$a_{ij},b_{ij}$ using the same map as in the proof of
Theorem~\ref{braidthm} in Section~\ref{ssec:braidproof}; then $\d$
and $\tilde{\d}$ agree on $\tilde{a}_{ij}$ and $\tilde{b}_{ij}$. Now
set
\begin{align*}
\tilde{C} &= \Phil{\sigma_k}(\hom{B}(A)) \cdot C \cdot \Phir{\sigma_k}(A)
+ \Theta^L_k(A) \cdot A \cdot \Phir{\sigma_k}(A), \\
\tilde{D} &= \Phil{\sigma_k}(A) \cdot D \cdot \Phir{\sigma_k}(\hom{B}(A))
+ \Phil{\sigma_k}(A) \cdot A \cdot \Theta^R_k(A), \\
\tilde{E} &= \Phil{\sigma_k}(\hom{B}(A)) \cdot E \cdot 
\Phir{\sigma_k}(\hom{B}(A))
- \Theta^L_k(A) \cdot D \cdot \Phir{\sigma_k}(\hom{B}(A))\\
&\quad\qquad + \Theta^L_k(A) \cdot \Theta^R_k(A), \\
\tilde{F} &= \Phil{\sigma_k}(\hom{B}(A)) \cdot F \cdot 
\Phir{\sigma_k}(\hom{B}(A))
+ \Phil{\sigma_k}(\hom{B}(A)) \cdot C \cdot \Theta^R_k(A)\\
&\quad\qquad + \Theta^L_k(A) \cdot (1+A) \cdot \Theta^R_k(A).
\end{align*}
We claim that these identifications make $\d$ and $\tilde{\d}$ agree
on $\tilde{c}_{ij},\tilde{d}_{ij},\tilde{e}_{ij},\tilde{f}_{ij}$,
which implies that $\d = \tilde{\d}$ on the entire algebra. We will
check that $\tilde{\d}\tilde{C} = \d \tilde{C}$ and
$\tilde{\d}\tilde{E} = \d\tilde{E}$; the proofs for $D$ and $F$ are
completely analogous.

By Proposition~\ref{chainrule} and Corollary~\ref{inverse}, we have
\begin{align*}
\Phil{\tilde{B}}(\tilde{A}) &=
\Phil{\sigma_k^{-1}B\sigma_k}(\hom{\sigma_k}(A))
= \Phil{B\sigma_k}(A) \cdot \Phil{\sigma_k^{-1}}(\hom{\sigma_k}(A))\\
&= \Phil{\sigma_k}(\hom{B}(A)) \cdot \Phil{B}(A) \cdot
(\Phil{\sigma_k}(A))^{-1},
\end{align*}
and similarly
\[
\Phir{\tilde{B}}(\tilde{A}) =
(\Phir{\sigma_k}(A))^{-1} \cdot \Phir{B}(A) \cdot
\Phir{\sigma_k}(\hom{B}(A)).
\]
Recall from Proposition~\ref{matrix} that
$\tilde{A} = \hom{\sigma_k}(A) = \Phil{\sigma_k}(A) \cdot A \cdot \Phir{\sigma_k}(A)$;
it follows that
\begin{align*}
\d \tilde{C} &= \Phil{\sigma_k}(\hom{B}(A)) \cdot
\left(1-\Phil{B}(A)\right)\cdot A \cdot \Phir{\sigma_k}(A)\\
&\qquad\qquad + \left(\Phil{\sigma_k}(A)-\Phil{\sigma_k}(\hom{B}(A))\right) 
\cdot A \cdot \Phir{\sigma_k}(A) \\
&= \left(\Phil{\sigma_k}(A) - \Phil{\sigma_k}(\hom{B}(A)) \cdot
\Phil{B}(A)\right)
\cdot A \cdot \Phir{\sigma_k}(A) \\
&= (1 - \Phil{\tilde{B}}(\tilde{A})) \cdot \tilde{A} \\
&= \tilde{\d}\tilde{C}
\end{align*}
and
\[
\begin{split}
\d\tilde{E} &= \Phil{\sigma_k}(\hom{B}(A)) \cdot (B - D - C \cdot
\Phir{B}(A)) \cdot \Phir{\sigma_k}(\hom{B}(A)) \\
& \qquad\qquad
- \left(\Phil{\sigma_k}(A)-\Phil{\sigma_k}(\hom{B}(A))\right) \cdot D \cdot \Phir{\sigma_k}(\hom{B}(A)) \\
& \qquad\qquad + \Theta^L_k(A) \cdot A \cdot (1 - \Phir{B}(A)) \cdot
\Phir{\sigma_k}(\hom{B}(A))
+ \d(\Theta^L_k(A) \cdot \Theta^R_k(A)) \\
&= \tilde{B} - \tilde{D} - \tilde{C} \cdot \Phir{\tilde{B}}(\tilde{A}) \\
&= \tilde{\d}\tilde{E},
\end{split}
\]
where we use the fact (a direct computation from the definition of
$\tilde{b}_{ij}$) that
\[
\begin{split}
\tilde{B} &= \Phil{\sigma_k}(\hom{B}(A)) \cdot B \cdot
\Phir{\sigma_k}(\hom{B}(A)) + \Phil{\sigma_k}(A) \cdot A \cdot
\Theta^R_k(A) \\
& \qquad + \Theta^L_k(A) \cdot A \cdot \Phir{\sigma_k}(\hom{B}(A)) +
\d\left( \Theta^L_k(A) \cdot \Theta^R_k(A) \right).
\end{split}
\]

Now that we have established that $\d = \tilde{\d}$ on the entire
algebra, it only remains to verify that the identification of
variables $A\mapsto\tilde{A}$, etc., constitutes a tame isomorphism.
This is straightforward; as in the proof of Theorem~\ref{braidthm},
the map $A\mapsto\tilde{A}, B\mapsto\tilde{B}$ is a tame
automorphism. We can then send $C$ to $\tilde{C}$ via the
composition of tame automorphisms
\begin{align*}
C &\mapsto C + \Theta^L_k(A) \cdot A \cdot \Phir{\sigma_k}(A), \\
C &\mapsto \Phil{\sigma_k}(\hom{B}(A)) \cdot C, \\
C &\mapsto C \cdot \Phir{\sigma_k}(A),
\end{align*}
of which the last two are tame because of Lemma~\ref{sigmak}, and
similarly for $D \mapsto \tilde{D}$. Finally, we can send $E$ to
$\tilde{E}$ via the composition of tame automorphisms
\begin{align*}
E &\mapsto E - \Theta^L_k(A) \cdot D \cdot \Phir{\sigma_k}(\hom{B}(A))
+ \Theta^L_k(A) \cdot \Theta^R_k(A), \\
E &\mapsto \Phil{\sigma_k}(\hom{B}(A)) \cdot E, \\
E &\mapsto E \cdot \Phir{\sigma_k}(\hom{B}(A)),
\end{align*}
and similarly for $F \mapsto \tilde{F}$.

\medskip
\textbf{Positive stabilization}\qua Let $B\in B_n$ be a
braid; as in the proof of Theorem~\ref{HCknotinvariance}, we denote
stabilization by adding a strand labelled $0$  and considering the
braid $\tilde{B} = B \sigma_0$. We will prove that the knot DGAs for
$B$ and $\tilde{B}$ are stable tame isomorphic.

Denote the knot DGA for $B$ by $(\A,\d)$ and the knot DGA for
$\tilde{B}$ by $(\tilde{\A},\tilde{\d})$, where the generators of
$\tilde{\A}$ have tildes for notational clarity. Thus $\tilde{\A}$
is generated by $\{\tilde{a}_{ij}\,|\,0\leq i\neq j\leq n\}$ in
degree $0$, and so forth. Note that $\tilde{\A}$ has more generators
than $\A$ does; we will establish that we can obtain
$(\tilde{\A},\tilde{\d})$ from $(\A,\d)$ by a suitable number of
stabilizations.

We will progressively identify generators of $\A$ and $\tilde{\A}$
so that $\d$ and $\tilde{\d}$ agree. The first identification is to
drop all tildes on the $\tilde{a}_{ij}$ generators; that is, set
$a_{ij} = \tilde{a}_{ij}$ for all $i,j\geq 0$, $i\neq j$. Note that
this introduces $2n$ generators $a_{0j}$ and $a_{i0}$ not in $\A$.

The matrices $\Phil{B\sigma_0}$ and $\Phir{B\sigma_0}$ were computed
in the proof of Theorem~\ref{HCknotinvariance}. We can then compute
$\tilde{\d}$ in terms of $\Phi^L = \Phil{B}$ and $\Phi^R =
\Phir{B}$. For $\tilde{b}_{ij}$ and $\tilde{c}_{ij}$, we already
calculated $\tilde{\d}$ in the proof of
Theorem~\ref{HCknotinvariance}. For completeness, we include here
the differentials of $\tilde{d}_{ij}$ and $\tilde{e}_{ij}$:
\begin{align*}
\tilde{\d}\tilde{d}_{00} &= \tilde{b}_{00} - \tilde{c}_{00} +
\tilde{b}_{00} \hom{B}(a_{01}) +
\tilde{b}_{0\ell} \Phi^R_{\ell 1} - \hom{B}(a_{10}) \tilde{c}_{00} - \Phi^L_{1\ell} \tilde{c}_{\ell 0} \\
\tilde{\d}\tilde{d}_{01} &= \tilde{b}_{01} - \tilde{c}_{01} -
\tilde{b}_{00} - \hom{B}(a_{10}) \tilde{c}_{01}
- \Phi^L_{1\ell} \tilde{c}_{\ell 1} \\
\tilde{\d}\tilde{d}_{0j} &= \tilde{b}_{0j} - \tilde{c}_{0j} -
\tilde{b}_{0\ell} \Phi^R_{\ell j}
- \hom{B}(a_{10}) \tilde{c}_{0j} - \Phi^L_{1\ell} \tilde{c}_{\ell j} \\
\tilde{\d}\tilde{d}_{10} &= \tilde{b}_{10} - \tilde{c}_{10} +
\tilde{c}_{00} + \tilde{b}_{10} \hom{B}(a_{01})
+ \tilde{b}_{1\ell} \Phi^R_{\ell 1} \\
\tilde{\d}\tilde{d}_{11} &= \tilde{b}_{11} - \tilde{c}_{11} - \tilde{b}_{10} + \tilde{c}_{01} \\
\tilde{\d}\tilde{d}_{1j} &= \tilde{b}_{1j} - \tilde{c}_{1j} + \tilde{c}_{0j} - \tilde{b}_{1\ell} \Phi^R_{\ell 1} \\
\tilde{\d}\tilde{d}_{j0} &= \tilde{b}_{j0} - \tilde{c}_{j0} + \tilde{b}_{j0} \hom{B}(a_{01})
+ \tilde{b}_{j\ell} \Phi^R_{\ell 1} + \Phi^L_{j\ell} \tilde{c}_{\ell 0} \\
\tilde{\d}\tilde{d}_{j1} &= \tilde{b}_{j1} - \tilde{c}_{j1} - \tilde{b}_{j0} + \Phi^L_{j\ell} \tilde{c}_{\ell 1} \\
\tilde{\d}\tilde{d}_{j_1j_2} &= \tilde{b}_{j_1j_2} - \tilde{c}_{j_1j_2} -
\tilde{b}_{j_1\ell} \Phi^R_{\ell j_2} + \Phi^L_{j_1\ell} \tilde{c}_{\ell j_2} \\
\tilde{\d}\tilde{e}_0 &= \tilde{b}_{00} - \hom{B}(a_{10}) \tilde{c}_{00} - \Phi^L_{1\ell} \tilde{c}_{\ell 0} \\
\tilde{\d}\tilde{e}_1 &= \tilde{b}_{11} + \tilde{c}_{01} \\
\tilde{\d}\tilde{e}_j &= \tilde{b}_{jj} + \Phi^L_{j\ell} \tilde{c}_{\ell j},
\end{align*}
where $j,j_1,j_2 \geq 2$ and, as usual, all monomials involving
$\ell$ are summed from $\ell=1$ to $\ell=n$.

If we set
\begin{xalignat*}{2}
b_{1i} &= - \tilde{b}_{0i} + \tilde{b}_{1i} - \tilde{b}_{00} a_{0i} & \qquad
c_{i1} &= - \tilde{c}_{i0} + \tilde{c}_{i1} - a_{i0} \tilde{c}_{00} \\
b_{ji} &= \tilde{b}_{ji} & \qquad
c_{ij} &= \tilde{c}_{ij}
\end{xalignat*}
for $i\geq 1$ and $j\geq 2$, then it is easy to check, from the
expressions for $\tilde{\d}\tilde{b}_{ij}$,
$\tilde{\d}\tilde{c}_{ij}$ given in the proof of
Theorem~\ref{HCknotinvariance}, that $\tilde{\d} b_{ij} = \d b_{ij}$
and $\tilde{\d} c_{ij} = \d c_{ij}$ for all $i,j\geq 1$. Next set
\begin{gather*}
\begin{split}
d_{11} &= \tilde{d}_{00} - \tilde{d}_{01} + \tilde{d}_{11} - \tilde{d}_{10} +
\tilde{b}_{00} \tilde{c}_{01} + \tilde{b}_{10} \tilde{c}_{00} \\
d_{1j} &= - \tilde{d}_{0j} + \tilde{d}_{1j} + \tilde{b}_{00} \tilde{c}_{0j} \\
d_{j1} &= - \tilde{d}_{j0} + \tilde{d}_{j1} + \tilde{b}_{j0} \tilde{c}_{00} \\
d_{j_1j_2} &= \tilde{d}_{j_1j_2} \\
e_1 &= \tilde{e}_0 + \tilde{e}_1 - \tilde{d}_{01} + \tilde{b}_{00} \tilde{c}_{01} \\
e_j &= \tilde{e}_j \\
d_{00} &= \tilde{d}_{10} - \tilde{d}_{11} + \tilde{d}_{01} - \tilde{b}_{10} \tilde{c}_{00}
- \tilde{b}_{00} \tilde{c}_{01} \\
e_0 &= \tilde{d}_{01} - \tilde{e}_1 - \tilde{b}_{00} \tilde{c}_{01}
\end{split} \\
d_{10} = \tilde{d}_{11} - \tilde{e}_1 \qquad
d_{01} = \tilde{e}_1 \qquad
d_{0j} = \tilde{d}_{1j} \qquad
d_{j0} = \tilde{d}_{j1}
\end{gather*}
for $j,j_1,j_2 \geq 2$, and then successively set
\begin{xalignat*}{3}
b_{10} &= \tilde{b}_{10} + \tilde{c}_{11} & \textrm{and} &&
c_{01} &= \tilde{c}_{01} + \tilde{b}_{11}; \\
b_{00} &= \tilde{b}_{00} - b_{11} - \Phi^L_{1\ell} \tilde{c}_{\ell 1} & \textrm{and} &&
c_{00} &= \tilde{c}_{00} - c_{11} - \tilde{b}_{1\ell} \Phi^R_{\ell 1}; \\
b_{j0} &= \tilde{b}_{j0} + \tilde{c}_{j1} - b_{j1} - \Phi^L_{j\ell} \tilde{c}_{\ell 1} & \textrm{and} &&
c_{0j} &= \tilde{c}_{0j} + \tilde{b}_{1j} - c_{1j} - \tilde{b}_{1\ell}\Phi^R_{\ell j}; \\
b_{0i} &= \tilde{b}_{1i} & \textrm{and} &&
c_{i0} &= \tilde{c}_{i1}
\end{xalignat*}
for $i\geq 1$ and $j \geq 2$.

Under the above identifications, the map sending the tilde variables
to the corresponding non-tilde variables is a tame automorphism on
$\tilde{\A}$. To see this, start with the generators of
$\tilde{\A}$, and now perform the following changes of variables.
First introduce new variables $b_{1i},c_{i1}$ for $i\geq 1$, and
then eliminate the generators $\tilde{b}_{0i},\tilde{c}_{i0}$ by
using the first identifications: $\tilde{b}_{0i} = - b_{1i} +
\tilde{b}_{1i} - \tilde{b}_{00} a_{0i}$ and $\tilde{c}_{i0} = -
c_{i1} + \tilde{c}_{i1} - a_{i0} \tilde{c}_{00}$. Now successively
continue this process using the given identifications, from first to
last; for each identification, introduce a new generator given by
the left hand side, and use the relation to eliminate the generator
given by the first term on the right hand side. This can be done
because no eliminated generator appears in successive
identifications. The final result is that all tilde variables are
eliminated and replaced by non-tilde variables. It follows that the
map sending tilde generators to non-tilde generators is tame.

On the other hand, a tedious but straightforward series of
calculations demonstrates that $\d = \tilde{\d}$ after the above
identifications have been made, where we set $\d b_{0i} = -a_{0i} +
a_{1i}$, $\d c_{i0} = -a_{i0} + a_{i1}$, $\d e_0 = b_{00}$, $\d
d_{00} = b_{00} - c_{00}$, $\d d_{i0} = -b_{i0}$, $\d d_{0i} =
c_{0i}$ for $i\geq 1$. It follows that $(\tilde{\A},\tilde{\d})$ is
tamely isomorphic to the result of stabilizing $(\A,\d)$ by
appending generators $a_{0i}$, $a_{i0}$, $b_{00}$, $b_{0i}$,
$b_{i0}$, $c_{00}$, $c_{0i}$, $c_{i0}$, $d_{00}$, $d_{0i}$,
$d_{i0}$, $e_0$, and extending $\d$ as in the previous sentence.
(This corresponds to $2n$ stabilizations of degree $0$ and $2n+2$
stabilizations of degree $1$.) This completes the proof that the
knot DGAs for $B$ and $\tilde{B}$ are stable tame isomorphic.

\medskip 
\textbf{Negative stabilization}\qua We now wish to prove that
the knot DGAs for $B$ and $\tilde{B} = B\sigma_0^{-1}$ are stable
tame isomorphic; this proof is entirely similar to the proof for
positive stabilization. We omit the calculations, but the sequence
of successive identifications is as follows, where indices satisfy
$i\geq 1$ and $j,j_1,j_2 \geq 2$, and all monomials involving $\ell$
are understood to be summed over $\ell\geq 1$:
\begin{xalignat*}{2}
b_{1i} &= \tilde{b}_{1i} - \tilde{b}_{0i} + \tilde{c}_{00}
\Phi^L_{1\ell} a_{\ell i} & \qquad
c_{i1} &= \tilde{c}_{i1} - \tilde{c}_{i0} + a_{i\ell} \Phi^R_{\ell 1} \tilde{b}_{00} \\
b_{ji} &= \tilde{b}_{ji} & \qquad
c_{ij} &= \tilde{c}_{ij},
\end{xalignat*}
followed by
\begin{gather*}
\begin{split}
d_{11} &= - \tilde{d}_{10} + \tilde{d}_{11} - \tilde{d}_{01} + \tilde{d}_{00}
- b_{1\ell} \Phi^R_{\ell 1} \tilde{b}_{00} - \tilde{c}_{00} \Phi^L_{1\ell} c_{\ell 1}
-2 \tilde{c}_{00}\tilde{b}_{00} \\
d_{1j} &= \tilde{d}_{1j} - \tilde{d}_{0j} - \tilde{c}_{00} \Phi^L_{1\ell} c_{\ell j} \\
d_{j1} &= \tilde{d}_{j1} - \tilde{d}_{j0} - b_{j\ell} \Phi^R_{\ell 1} \tilde{b}_{00} \\
d_{j_1j_2} &= \tilde{d}_{j_1j_2} \\
e_1 &= -\tilde{d}_{01} + \tilde{e}_0 + \tilde{e}_1 - \tilde{c}_{00} \Phi^L_{1\ell} c_{\ell 1}
- \tilde{c}_{00} \tilde{b}_{00} + (\tilde{e}_0-\tilde{d}_{00})\hom{B}(a_{10})
+ \hom{B}(a_{01})\tilde{e}_0 \\
e_j &= \tilde{e}_j \\
\end{split} \\
d_{00} = \tilde{d}_{00} \quad\qua
e_0 = \tilde{e}_0 \quad\qua
d_{10} = \tilde{d}_{11} - \tilde{e}_1 \quad\qua
d_{01} = \tilde{e}_1 \quad\qua
d_{j0} = \tilde{d}_{j0} \quad\qua
d_{0j} = \tilde{d}_{0j},
\end{gather*}
and finally
\begin{xalignat*}{3}
\begin{split}
b_{10} &= -\tilde{b}_{10} + c_{11} + \tilde{c}_{10} - a_{1\ell} \Phi^R_{\ell 1} \tilde{b}_{00} \\
  -& (b_{1\ell}\Phi^R_{\ell 1} + \tilde{b}_{0\ell}\Phi^R_{\ell
1} + 2 \tilde{c}_{00}) \hom{B}(a_{10})
\end{split}
& \textrm{and} &&
\begin{split}
c_{01} &= -\tilde{c}_{01} + b_{11} + \tilde{b}_{01} - \tilde{c}_{00}\Phi^L_{1\ell} a_{\ell 1} \\
 - &\hom{B}(a_{01}) (\Phi^L_{1\ell} c_{\ell 1} +
\Phi^L_{1\ell} \tilde{c}_{\ell 0} + 2 \tilde{b}_{00});
\end{split} \\
b_{00} &= \tilde{b}_{00} + \Phi^L_{1\ell} \tilde{c}_{\ell 0} & \textrm{and} &&
c_{00} &= \tilde{c}_{00} + \tilde{b}_{0\ell} \Phi^R_{\ell 1}; \\
%b_{10} &= \tilde{b}_{10} - c_{11} - \tilde{c}_{10}
%+ \tilde{b}_{0\ell} \Phi^R_{\ell 1} \hom{B}(a_{10}) & \textrm{and} &&
%c_{01} &= \tilde{c}_{01} - b_{11} - \tilde{b}_{01} -
%\hom{B}(a_{01}) \Phi^L_{1\ell} \tilde{c}_{\ell 0}; \\
b_{j0} &= - \tilde{b}_{j0} + \tilde{c}_{j0}
+ b_{j\ell} \Phi^R_{\ell 1} - \Phi^L_{j\ell} \tilde{c}_{\ell 0} & \textrm{and} &&
c_{0j} &= - \tilde{c}_{0j} + \tilde{b}_{0j}
+ \Phi^L_{1\ell} c_{\ell j} - \tilde{b}_{0\ell} \Phi^R_{\ell j}; \\
b_{0i} &= \tilde{b}_{0i} & \textrm{and} &&
c_{i0} &= \tilde{c}_{i0}.
\end{xalignat*}
We find that $\d=\tilde{\d}$ as before, where we extend $\d$ by defining
$\d b_{0i} = a_{0i} - \Phi^L_{1\ell} a_{\ell i}$,
$\d c_{i0} = a_{i0} - a_{i\ell} \Phi^R_{\ell 1}$,
$\d d_{00} = b_{00}-c_{00}$, $\d e_0 = b_{00}$,
$\d d_{i0} = - b_{i0}$, $\d d_{0i} = c_{0i}$.
(The calculation uses the fact that
$\Phi^L_{1\ell} a_{\ell m} \Phi^R_{m1} = (\hom{B}(A))_{11} = -2$.)
It follows that the knot DGA for $B \sigma_0^{-1}$ is tamely isomorphic
to a stabilization of the knot DGA for $B$, as desired. \qed

%*********************************************************************
%*********************************************************************
\section{Properties of the invariants}
\label{sec:properties}

%*********************************************************************
\subsection{Conjugation and abelianization}
\label{ssec:conjugation}

There is a symmetry in the homomorphism $\phi\co  B_n\to\Aut(\A_n)$
which induces a symmetry in braid and knot DGAs. Define an
involution on $\A_n$, which we term \textit{conjugation} and write
as $v\mapsto\overline{v}$, as follows: the conjugate of a monomial
$a_{i_1j_1}\cdots a_{i_mj_m}$ is $a_{j_mi_m}\cdots a_{j_1i_1}$;
extend conjugation linearly, and set $\overline{1}=1$. In other
words, conjugation replaces each generator $a_{ij}$ by $a_{ji}$ and
then reverses the order of each word. The key observation is the
following.

\begin{proposition}
For all $v\in\A_n$ and any $B\in B_n$, we have
$\hom{B}(\overline{v}) = \overline{\hom{B}(v)}$.
\label{conjprop}
\end{proposition}

\begin{proof}
The result holds for $B=\sigma_k$ and $v=a_{ij}$ for any $i,j,k$, by
direct inspection of the definition of $\phi$. The general result
for $B=\sigma_k$ follows since $\overline{vw} = \overline{w} \,
\overline{v}$ for all $v,w\in\A_n$; since $\phi$ is a homomorphism,
the entire proposition follows.
\end{proof}

One consequence of Proposition~\ref{conjprop} is that the matrices
$\Phil{B}(A),\Phir{B}(A)$ determine each other. Extend the operation
of conjugation to matrices by conjugating each entry, and let
$^{\transpose}$ denote transpose; then we have the following result.

\begin{proposition}
For any $B\in B_n$, $\Phir{B}(A) =
\overline{\Phil{B}(A)}^{\transpose}$.
\label{philphir}
\end{proposition}

\begin{proof}
By Proposition~\ref{conjprop}, we have
\[
\hom{B}(a_{\dott i}) = \hom{B}(\overline{a_{i\dott}}) = \overline{\hom{B}(a_{i\dott})}
= \overline{ \sum_{j=1}^n (\Phil{B}(A))_{ij}a_{j\dott} }
= \sum_{j=1}^n  a_{\dott j} (\overline{\Phil{B}(A)})_{ij};
\]
now use the definition of $\Phir{B}(A)$.
\end{proof}

Conjugation also allows us to define simplified commutative versions
of the braid and knot DGAs. Intuitively, we can quotient the braid
and knot DGAs by conjugation and abelianize; the differential is
well-defined on the resulting quotient algebras. For instance, on
the braid DGA, if $\d b_{ij} = a_{ij} - \hom{B}(a_{ij})$, then $\d
b_{ji} = a_{ji} - \hom{B}(a_{ji}) =
\overline{a_{ij}-\hom{B}(a_{ij})}$, and so $\d$ is still
well-defined if we mod out by $a_{ji} = a_{ij}$, $b_{ji} = b_{ij}$
and then abelianize.

\begin{definition}
Let $B\in B_n$ be a braid. Let $\A$ be the graded sign-commutative
algebra on $n(n-1)$ generators, $\{a_{ij}\,|\,1\leq i<j\leq n\}$ of
degree $0$ and $\{b_{ij}\,|\,1\leq i<j\leq n\}$ of degree $1$.
Define a differential $\d$ on $\A$ to be the usual differential on
the braid DGA, where we set $a_{ji} = a_{ij}$ for $j>i$. Then we
call $(\A,\d)$ the \textit{abelian braid DGA} of $B$.
\end{definition}

\noindent Here by ``sign-commutative,'' we mean that $vw =
(-1)^{(\deg v)(\deg w)}wv$ for all $v,w\in\A$ of pure degree.

For the knot DGA, extend the conjugation operation to the entire
algebra by defining $\overline{b_{ij}} = c_{ji}$, $\overline{c_{ij}}
= b_{ji}$, $\overline{d_{ij}} = -d_{ji}$, $\overline{e_i} = e_i -
d_{ii}$, and $\overline{vw} = (-1)^{(\deg v)(\deg w)}
\overline{w}\,\overline{v}$.

\begin{proposition}
If $(\A,\d)$ is the knot DGA for some braid, then $\overline{\d v} = \d \overline{v}$ for
all $v\in\A$.
\end{proposition}

\begin{proof}
We calculate directly that
\begin{gather*}
\overline{\d b_{ij}} = \overline{a_{ij}} - \sum_\ell
\overline{a_{\ell j}}\, \overline{(\Phil{B}(A))_{i\ell}}
= \d c_{ji}, \\
\overline{\d c_{ij}} = \overline{a_{ij}} - \sum_\ell
\overline{(\Phir{B}(A))_{\ell j}} \, \overline{a_{i\ell}}
= \d b_{ji}, \\
\overline{\d d_{ij}} = \overline{b_{ij}}-\overline{c_{ij}}+\sum_\ell
\left( - \overline{(\Phir{B}(A))_{\ell j}} \, \overline{b_{i\ell}}
+ \overline{c_{\ell j}} \, \overline{(\Phil{B}(A))_{i\ell}} \right) %\\
= -\d d_{ji}, \\
\overline{\d e_i} = \overline{b_{ii}}+\sum_\ell \overline{c_{\ell
i}} \, \overline{(\Phil{B}(A))_{i\ell}} = c_{ii} + b_{i\ell}
(\Phir{B}(A))_{\ell i} = \d e_i - \d d_{ii},
\end{gather*}
and if $\overline{\d v} = \d \overline{v}$ and $\overline{\d w} = \d
\overline{w}$, then
\[
\overline{\d (vw)} = (-1)^{(\deg v - 1)(\deg w)} \overline{w} \,
\overline{\d v} + (-1)^{(\deg v)(\deg w)} \overline{\d w} \,
\overline{v} = \d (\overline{vw}).
\]
The proposition follows.
\end{proof}

Since $\d$ commutes with conjugation, we may mod out by conjugation
and abelianize, and $\d$ is well-defined on the resulting algebra.

\begin{definition}
Let $B\in B_n$. Let $\A$ be the graded sign-commutative algebra with
generators $\{a_{ij}\,|\,1\leq i<j\leq n\}$ of degree $0$,
$\{b_{ij}\,|\,1\leq i,j\leq n\}$ of degree $1$, and
$\{d_{ij}\,|\,1\leq i<j\leq n\}$ and $\{e_i\,|\,1\leq i\leq n\}$ of
degree $2$.
\label{abeliandef}
Define $\d$ on $\A$ to be the usual differential on the knot DGA,
where we set $a_{ji} = a_{ij}$ for $j>i$ and $c_{ij} = b_{ji}$ for
all $i,j$. Then we call $(\A,\d)$ the \textit{abelian knot DGA} of
$B$.
\end{definition}

The notion of stable tame isomorphism can be defined for polynomial
algebras just as for tensor algebras, and with respect to this
relation, abelian braid and knot DGAs give invariants.

\begin{proposition}
Up to stable tame isomorphism, the abelian braid DGA of $B$ is an
invariant of the conjugacy class of $B$, and the abelian knot DGA of
$B$ is an invariant of the knot closure of $B$.
\end{proposition}

\begin{proof}
Trace through the invariance proofs for the original knot and braid
DGAs, and note that all of the tame isomorphisms and stabilizations
used there are invariant under conjugation and hence well-defined in
the commutative case.
\end{proof}

As in the noncommutative case, our main interest in the abelian DGAs
lies in the homology in degree $0$, which we write as
$HC_0^{\textrm{ab}}(B)$ in the braid case and
$HC_0^{\textrm{ab}}(K)$ in the knot case.

\begin{corollary}
Up to isomorphism, $HC_0^{\operatorname{ab}}(B)$ is an invariant of
the conjugacy class of $B$, and $HC_0^{\operatorname{ab}}(K)$ is an
invariant of the knot $K$.
\end{corollary}

We can summarize the various flavors of $HC_0$ in a diagram, as
follows. Suppose that braid $B$ has closure $K$, and let
$\ab(HC_0(B)),\ab(HC_0(K))$ represent the abelianizations of
$HC_0(B),HC_0(K)$. Then $HC_0(B),HC_0(K)$ can both be written as
$\Z\langle \{a_{ij}\}\rangle$ modulo some ideal, and the ideal for
$HC_0(K)$ contains the ideal for $HC_0(B)$ by
Proposition~\ref{newdga}. Also, $\ab(HC_0(B)),HC_0^{\ab}(B)$ are
both $\Z[\{a_{ij}\}]$ modulo some ideal, and the ideal for
$HC_0^{\ab}(B)$ is generated by the ideal for $\ab(HC_0(B))$ and the
relations $\{a_{ij}-a_{ji}\}$. Thus we have the following
commutative diagram of invariants:
\[
\xymatrix{
HC_0(B) \ar@{->>}[r] \ar@{->>}[d] & \ab(HC_0(B)) \ar@{->>}[r] \ar@{->>}[d] &
HC_0^{\ab}(B) \ar@{->>}[d] \\
HC_0(K) \ar@{->>}[r] & \ab(HC_0(K)) \ar@{->>}[r] & HC_0^{\ab}(K).
}
\]

%*********************************************************************
\subsection{Mirrors and inverses}
\label{ssec:mirrors}

In this section, we show that the braid and knot DGAs do not change
under mirroring and inversion. Recall that the \textit{mirror} of a
knot $K$ is the knot $\overline{K}$ obtained by reversing each of
the crossings of the diagram of $K$, and similarly for the mirror
$\overline{B}$ of a braid $B$. In $B_n$, mirroring is the
homomorphism sending $\sigma_k$ to $\sigma_k^{-1}$ for each $k$. It
will be convenient for us to consider a related mirror operation on
$B_n$, $B \mapsto B^*$, which is the homomorphism sending $\sigma_k$
to $\sigma_{n-k}^{-1}$ for each $k$. Recall also that the
\textit{inverse} of a knot is the knot with the opposite
orientation. If $K$ is the closure of $B\in B_n$, then the mirror
$\overline{K}$ is the closure of both $\overline{B}$ and $B^*$,
while the inverse of $K$ is the closure of $\overline{B}^{-1}$.

\begin{proposition}
The braid DGAs for a braid $B$, its mirrors $\overline{B}$ and
$B^*$, and its inverse $B^{-1}$ are stable tame isomorphic.
\label{braidmirror}
\end{proposition}

\begin{proof}
Let $B\in B_n$. We first show that $B$ and $B^{-1}$ have equivalent
braid DGAs. In the braid DGA of $B^{-1}$, we have $\d b_{ij} =
a_{ij}-\hom{B^{-1}}(a_{ij})$ for all $i,j$; hence if we replace
$a_{ij}$ by $\hom{B}(a_{ij})$ and $b_{ij}$ by $-b_{ij}$ for all
$i,j$, then we obtain $\d b_{ij} = a_{ij}-\hom{B}(a_{ij})$, which
gives the braid DGA for $B$.

Next, as in \cite{Bir}, denote
$\Delta=(\sigma_1\cdots\sigma_{n-1})(\sigma_1\cdots\sigma_{n-2})\cdots(\sigma_1\sigma_2)\sigma_1$.
Then $\sigma_{n-k} = \Delta\sigma_k\Delta^{-1}$ for all $k$, and so
$\overline{B}$ and $B^*$ are conjugate for any $B\in B_n$; hence by
Theorem~\ref{braidthm}, the braid DGAs of $\overline{B}$ and $B^*$
are stable tame isomorphic.

It now suffices to show that $B$ and $B^*$ have equivalent braid
DGAs. Consider the tame automorphism $\xi$ of the algebra $\A_n$
which sends $a_{ij}$ to $a_{n+1-i,n+1-j}$ for each generator
$a_{ij}$ (and sends $1$ to $1$). By inspecting the definition of
$\phi$, we see that the induced action of $\xi$ sends
$\phi_{\sigma_k}$ to $\phi_{\sigma_{n-k}^{-1}}$, ie, $\xi
\phi_{\sigma_k} = \phi_{\sigma_{n-k}^{-1}} \xi$ for all $k$. Hence
$\xi \phi_{B^*} = \phi_B \xi$ for all $B$.

Write the braid DGAs for $B,B^*$ as $(\A,\d),(\A,\tilde{\d})$,
respectively, and extend $\xi$ to a tame automorphism of $\A$ by
sending $b_{ij}$ to $b_{n+1-i,n+1-j}$ for each $b_{ij}$. Then
\[
\d \xi b_{ij} = \d b_{n+1-i,n+1-j} = \xi(a_{ij}) - \phi_B(\xi(a_{ij})) =
\xi(a_{ij} - \phi_{B^*}(a_{ij})) = \xi \tilde{\d} b_{ij},
\]
and so $\xi$ intertwines the differentials $\d$ and $\tilde{\d}$.
\end{proof}

\begin{proposition}
The knot DGA classes for a knot, its mirror, and its inverse are the same.
\label{knotmirror}
\end{proposition}

\begin{proof}
Let $B \in B_n$ be a braid whose closure is a knot $K$. We first
show that the knot DGAs for $B,B^*$ are stable tame isomorphic,
which implies that $K$ and its mirror have the same knot DGA class.
Denote these knot DGAs by $(\A,\d),(\A,\d')$, respectively.

Let $\Xi$ be the map on $n\times n$ matrices which permutes entries
according to the definition $(\Xi(M))_{ij} = M_{n+1-i,n+1-j}$; note
that $\Xi$ commutes with matrix multiplication. Define the tame
automorphism $\xi$ on $\A$ by its action on generators: $\xi(A) =
\Xi(A)$, $\xi(B) = \Xi(B)$, $\xi(C) = \Xi(C)$, $\xi(D) = \Xi(D)$,
$\xi(e_i) = e_{n+1-i}$.

Inspection of the definitions of $\phi$ and $\Phi^L$ shows that
$\xi(\Phil{\sigma_k}(A)) = \Xi(\Phil{\sigma_{n-k}^{-1}}(A))$ for all
$k$. This, along with Proposition~\ref{chainrule} and the identity
$\xi \phi_{B^*} = \phi_B \xi$ from the proof of
Proposition~\ref{braidmirror}, allows us to prove readily by
induction that $\xi(\Phil{B^*}(A)) = \Xi(\Phil{B}(A))$. Similarly,
we can establish that $\xi(\Phir{B^*}(A)) = \Xi(\Phir{B}(A))$.

Hence
\[
\d\xi B = \Xi\d B %= \Xi((1-\Phil{B}(A))\cdot A)
= (1-\Xi(\Phil{B}(A))) \cdot \Xi(A)
= (1-\xi(\Phil{B^*}(A))) \cdot \xi(A)
= \xi\d'B
\]
and similarly $\d\xi C  = \xi\d'C$, while
\[
\d\xi D = \Xi\d D
%= \Xi(B\cdot (1-\Phir{B}(A)) - (1-\Phil{B}(A))\cdot C)
= \xi(B) \cdot (1-\xi(\Phir{B^*}(A))) -
(1-\xi(\Phil{B^*}(A)))\cdot \xi(C) = \xi\d'D
\]
and
\[
\d \xi e_i = \d e_{n+1-i} = \left(\Xi(B+\Phil{B}(A)\cdot C)\right)_{ii}
= \left(\xi(B) + \xi(\Phil{B^*}(A))\cdot \xi(C)\right)_{ii}
= \xi \d'e_i.
\]
It follows that $\xi$ intertwines $\d$ and $\d'$, and so the knot DGAs of $B$ and $B^*$ are
equivalent.

Now we show that the knot DGAs of $B$ and $B^{-1}$ are equivalent;
this implies that $K$ and its inverse have the same knot DGA class,
since $B^{-1},(B^{-1})^*$ have equivalent knot DGAs by the above
argument, and $(B^{-1})^*,\overline{B}^{-1}$ are conjugate and hence
have equivalent knot DGAs by Theorem~\ref{knotthm}. Write the knot
DGA of $B^{-1}$ as $(\tilde{\A},\tilde{\d})$, where the generators
of $\tilde{\A}$ are marked by tildes. Identify the generators of
$\tilde{\A}$ with those of $\A$ by $\tilde{A}=\hom{B}(A)$,
$\tilde{B}=-B\cdot\Phir{B}(A)$, $\tilde{C}=-\Phil{B}(A)\cdot C$,
$\tilde{D}=D$, and $\tilde{e}_i=d_{ii}-e_i$. Using
Corollary~\ref{inverse} and Proposition~\ref{matrix}, we calculate
that
\[
\tilde{\d}\tilde{B} = (1-\Phil{B^{-1}}(\tilde{A}))\cdot\tilde{A} =
\left(1-(\Phil{B}(A))^{-1}\right)\cdot \Phil{B}(A) \cdot A \cdot
\Phir{B}(A) = \d\tilde{B}
\]
and similarly $\tilde{\d}\tilde{C}=\d\tilde{C}$, while
\[
\tilde{\d}\tilde{D} = \tilde{B}\cdot (1-\Phir{B^{-1}}(\tilde{A})) -
(1-\Phil{B^{-1}}(\tilde{A}))\cdot \tilde{C} =
B\cdot (1-\Phir{B}(A))-(1-\Phil{B}(A))\cdot C = \d\tilde{D}
\]
and
\[
\tilde{\d}\tilde{e}_i = (\tilde{B}+\Phil{B^{-1}}(\tilde{A})\cdot\tilde{C})_{ii}
= (-B\cdot\Phir{B}(A)-C)_{ii} = \d(d_{ii}-e_i) = \d\tilde{e}_i.
\]
Hence our identification of generators, which yields a tame
automorphism between $\A$ and $\tilde{\A}$, gives $\d=\tilde{\d}$,
and so the knot DGAs of $B$ and $B^{-1}$ are equivalent.
\end{proof}

It can similarly be shown that the abelian braid and knot DGAs do
not distinguish between mirrors or inverses. We may then ask how our
invariants relate to classical invariants which also fail to
distinguish between mirrors or inverses, such as the Alexander
polynomial and the absolute value of the signature. In
Section~\ref{ssec:det}, we establish a connection between the knot
DGA and the Alexander polynomial. First, in the next section, we
deduce an invariant from the DGAs which shows that our invariants
can distinguish between knots with the same Alexander polynomial and
signature.

%*********************************************************************
\subsection{Augmentation number}
\label{ssec:augment}

We now introduce a family of easily computable invariants derived
from the braid and knot DGAs. These invariants, which we call
augmentation numbers, gives a crude but effective way to distinguish
between braid conjugacy classes or knots. The advantage of
augmentation numbers is that, unlike $HC_0$ or the full DGA, they
are positive integers and hence require no ad hoc methods to
interpret (cf\ the computations of $HC_0$ in
Section~\ref{sec:computations}).

\begin{definition}
If $R$ is a finitely generated ring and $d\geq 2$ an integer, define
the $\Z_d$ \textit{augmentation number} of $R$ to be the number of
ring homomorphisms $R \to \Z_d$. For $B$ any braid and $K$ any knot,
write $\Aug(B,d)$, $\Augab(B,d)$, $\Aug(K,d)$, $\Augab(K,d)$ for the
$\Z_d$ augmentation numbers of $HC_0(B)$, $HC_0^{\ab}(B)$,
$HC_0(K)$, $HC_0^{\ab}(K)$, respectively.
\end{definition}

We use the term ``augmentation'' in keeping with \cite{EFM}, in
which the case $d=2$ is discussed in the context of DGAs and
Poincar\'e polynomials. One could imagine replacing $\Z_d$ by any
finite ring, although the usefulness of the more general
construction is unclear.

Calculating augmentation numbers is straightforward. For instance,
to compute $\Aug(B,d)$ for $B\in B_n$, write
$HC_0(B)=\Z\langle\{a_{ij}\}\rangle/\I$ and count the number of
points in $\Z_d^{n(n-1)}$ (the space parametrizing ring
homomorphisms $\Z\langle\{a_{ij}\}\rangle \to \Z_d$) which lie in
the vanishing set of $\I$. There are $d^{n(n-1)}$ possible points
and $n(n-1)$ conditions to check, corresponding to $\d b_{ij}$ for
$1\leq i\neq j\leq n$. Similarly, for $\Augab(B,d)$, there are
$d^{n(n-1)/2}$ possible points and $n(n-1)/2$ conditions; for
$\Aug(K,d)$ with $K$ the closure of $B$, $d^{n(n-1)}$ points and
$2n^2$ conditions; and for $\Augab(K,d)$, $d^{n(n-1)/2}$ points and
$n^2$ conditions. (Note thus that $\Augab$ typically takes
considerably less computing time than $\Aug$ to calculate.) In all
cases, the $\Z_d$ augmentation number is at least $1$, as the
following result shows.

\begin{proposition}
For any braid $B$ or knot $K$, the ring homomorphism from $\A_n$ to
$\Z$ defined by sending each generator $a_{ij}$ to $-2$ descends to
homomorphisms from $HC_0(B)$, $HC_0^{\ab}(B)$, $HC_0(K)$,
$HC_0^{\ab}(K)$ to $\Z$.
\label{prop:augment}
\end{proposition}

\begin{proof}
It suffices to prove the assertion for $HC_0^{\ab}(K)$, where $K$ is
the closure of the braid $B$; see the diagram at the end of
Section~\ref{ssec:conjugation}. Let $M$ denote the $n\times n$
matrix all of whose entries are $-2$; then we need to show that
$M=\Phil{B}(M)\cdot M$. It is clear from the definition of $\phi$
and Lemma~\ref{sigmak} that $\hom{\sigma_k}(M) = M$ and
$\Phil{\sigma_k}(M)\cdot M = M$. By Proposition~\ref{chainrule}, it
follows that
\[
\Phil{\sigma_k B}(M)\cdot M = \Phil{B}(\hom{\sigma_k}(M)) \cdot \Phil{\sigma_k}(M) \cdot M
= \Phil{B}(M) \cdot M
\]
for any $B$, and hence $\Phil{B}(M)\cdot M = M$ for all $B$.
\end{proof}

As an example of the computation of augmentation numbers, consider
the trefoil $3_1$, for which we calculated in
Section~\ref{ssec:HCdef} that $HC_0\cong \Z[x]/((x+2)(x-1))$; it is
also easy to see that $HC_0^{\ab}(3_1) \cong \Z[x]/((x+2)(x-1))$.
Then $\Aug(3_1,d)=\Augab(3_1,d)=2$ for $d \neq 3$ and
$\Aug(3_1,3)=\Augab(3_1,3)=1$. Note that these augmentation numbers
could be calculated without a nice form for $HC_0$, by directly
using the definition of $\d$.

Any knot $K$ with $HC_0(K) \cong \Z[x]/(p(x))$ for some polynomial
$p(x)$ satisfies $\Augab(K,d) \leq \Aug(K,d) \leq d$ for all $d$.
Hence augmentation numbers can be used to show that some knots have
$HC_0$ which is not of the form $\Z[x]/(p(x))$; see the examples in
Sections~\ref{ssec:connectsum} and~\ref{ssec:HOMFLY}.

%*********************************************************************
%*********************************************************************
\section{Knot contact homology and the determinant}
\label{sec:determinant}

For knots $K$, contact homology is closely related to the classical
invariant $|\Delta_K(-1)|$, the determinant of $K$. Throughout this
section, we will denote by $\Sigma_2 = \Sigma_2(K)$ the double
branched cover of $S^3$ over $K$, as in \cite{Rol}; it is well known
that $|H_1(\Sigma_2(K))|=|\Delta_K(-1)|$. In Section~\ref{ssec:det},
we establish a surjection from $HC_0(K)$ to a quotient of the
polynomial ring $\Z[x]$ specified by the largest invariant factor of
$H_1(\Sigma_2(K))$, which typically is equal to the determinant of
$K$. In Section~\ref{ssec:lin}, we discuss a natural linearization
of knot contact homology, and see that $H_1(\Sigma_2(K))$, and hence
the determinant, can be deduced from this linearization.

%*********************************************************************
\subsection{$HC_0(K)$ and $\Sigma_2(K)$}
\label{ssec:det}

Before stating the main result of this section, we need to introduce a sequence of polynomials
which plays a key role in the remainder of the paper. Define $\{p_m\in\Z[x]\}$ inductively by
$p_0(x)=2-x$, $p_1(x)=x-2$,
$p_{m+1}(x) = xp_m(x) - p_{m-1}(x)$.
Perhaps a more illuminating way to view $p_m$, $m \geq 1$, is via its roots:
\[
p_m(x) = \prod_{k=0}^{m-1} \left( x - 2\cos\frac{2k\pi}{2m-1} \right).
\]

\begin{theorem}
For $K$ a knot, let $n(K)$ denote the largest invariant factor of
the abelian group $H_1(\Sigma_2(K))$. Then there is a surjective
ring homomorphism
\[
HC_0(K) \twoheadrightarrow \Z[x]/(p_{(n(K)+1)/2}(x)),
\]
and this map factors through $HC_0^{\ab}(K)$. In particular, the
ranks of $HC_0(K)$ and $HC_0^{\ab}(K)$ are at least $(n(K)+1)/2$.
\label{detthm}
\end{theorem}

\noindent Note that $n(K)$ is odd since the determinant of a knot is
always odd, and that $n(K)$ divides $|\Delta_K(-1)|$ and is
divisible by any prime dividing $|\Delta_K(-1)|$. A table of
invariant factors of $H_1(\Sigma_2(K))$ for small knots can be found
in \cite{BZ}.

To prove Theorem~\ref{detthm}, we will need a number of preliminary
results. Define another sequence $\{q_m \in \Z[x]\}$ by $q_0(x)=-2$,
$q_1(x)=-x$, $q_{m+1}(x) = xq_m(x) - q_{m-1}(x)$. Note that the
recursion can be used to define $q_m$ for all $m\in\Z$, and that
$q_{-m} = q_m$ for all $m$ because of symmetry.

\begin{lemma}
For all $m_1,m_2\in\Z$, we have \label{funceq}
\[
q_{m_1}(x) q_{m_2}(x) + q_{m_1+m_2}(x) + q_{m_1-m_2}(x) = 0.
\]
\end{lemma}

\begin{proof}
The identity holds by the defining relations for $q_m$ for $m_2=0$
and $m_2=1$. It is then easy to establish by induction that it holds
for all $m_2 \geq 0$, and hence for all $m_2$ since $q_{-m_2} =
q_{m_2}$.
\end{proof}

\begin{lemma}
\label{divisible}
For all $\ell,m\in\Z$ with $m\geq 0$, we have $p_m(x) \, | \,
q_{\ell+2m-1}(x) - q_{\ell}(x)$.
\end{lemma}

\begin{proof}
An easy induction shows that $(p_m(x))^2 = (x+2)(q_{2m-1}(x)+2)$ and
$p_m(x)p_{m+1}(x) = -(x+2)(q_{2m}(x)-x)$ for all $m \geq 0$; since
$x+2$ divides $p_m(x),p_{m+1}(x)$, this establishes the lemma for
$\ell=0,1$. The lemma now follows by induction on $\ell$.
\end{proof}

One can more elegantly prove Lemmas~\ref{funceq} and~\ref{divisible}
by noting that the polynomials $p_m,q_m$ are characterized by the
identities $p_m(x+x^{-1}) = x^m+x^{-m}-x^{m-1}-x^{1-m}$,
$q_m(x+x^{-1}) = -x^m-x^{-m}$; the lemmas then follow easily.

We next need to recall the Burau representation on $B_n$; see
\cite{Bir} as a reference. For $1\leq k\leq n-1$, define
$\Burau{\sigma_k}$ to be the linear map on the free $\Z[t,t^{-1}]$
module on $n$ generators whose matrix is the identity except for the
$2\times 2$ submatrix formed by the $k,k+1$ rows and columns, which
is
$\left( \begin{smallmatrix} 1-t & t \\ 1 & 0 \end{smallmatrix} \right)$.
This extends to a representation of $B_n$, classically known as the
(nonreduced) \textit{Burau representation}.  The Burau
representation is reducible to a trivial $1$-dimensional
representation spanned by the vector $(1,\dots,1)$, and the
$(n-1)$-dimensional \textit{reduced Burau representation}. More
concretely, the map $\Z[t,t^{-1}]^n \to \Z[t,t^{-1}]^{n-1}$ given by
$(v_1,\dots,v_n) \mapsto (v_2-v_1,\dots,v_n-v_1)$ induces the
reduced Burau representation on $\Z[t,t^{-1}]^{n-1}$.

For the remainder of this section, we will set $t=-1$; then the
nonreduced and reduced Burau representations of $B_n$ act on $\Z^n$
and $\Z^{n-1}$, respectively. Let $\Bur{B}$ and $\red{B}$ denote the
matrices for $B$ in the two respective representations. It is well
known that $\red{B}$ is a presentation matrix for the group
$H_1(\Sigma_2(K))$.

We need a bit more notation. If $B$ is a braid, let $\hat{B}$ denote
the ``reverse'' braid obtained by reversing the word which gives
$B$; that is, if $B = \sigma_{i_1}^{j_1} \sigma_{i_2}^{j_2} \cdots
\sigma_{i_m}^{j_m}$, then $\hat{B} = \sigma_{i_m}^{j_m} \cdots
\sigma_{i_2}^{j_2} \sigma_{i_1}^{j_1}$. (In our notation from
Section~\ref{ssec:mirrors}, $\hat{B}=\overline{B}^{-1}$.) Also, for
$v,w \in \Z^n$, let $M(v,w)$ denote the $n\times n$ matrix defined
by $(M(v,w))_{ij} = q_{v_i - w_j}$.

\begin{lemma}
Let $B\in B_n$ and $v\in\Z^n$.  If we set $A = M(v,v)$, then \label{qlemma1}
$\hom{B}(A) = M(\Burau{\hat{B}}v,\Burau{\hat{B}}v)$.
\end{lemma}

\begin{proof}
Since $A = M(v,v)$, we have $a_{ij} = q_{v_i-v_j}$ for all $i,j$.
The key is to verify the
relation for $B = \sigma_k$. If $i\neq k,k+1$, then
\begin{align*}
\hom{\sigma_k}(a_{ki}) &= -a_{k+1,i} - a_{k+1,k} a_{ki} =
-q_{v_{k+1}-v_i} - q_{v_{k+1}-v_k} q_{v_k-v_i} =
q_{2v_k-v_{k+1}-v_i} \\
&= (M(\Burau{\sigma_k}v,\Burau{\sigma_k}v))_{ki}
\end{align*}
by Lemma~\ref{funceq}, and similarly $\hom{\sigma_k}(a_{ik}) =
(M(\Burau{\sigma_k}v,\Burau{\sigma_k}v))_{ik}$;
\[
\hom{\sigma_k}(a_{k+1,i}) = a_{ki} = q_{v_k-v_i} = (M(\Burau{\sigma_k}v,\Burau{\sigma_k}v))_{k+1,i}
\]
and similarly $\hom{\sigma_k}(a_{i,k+1}) = (M(\Burau{\sigma_k}v,\Burau{\sigma_k}v))_{i,k+1}$;
\[
\hom{\sigma_k}(a_{k,k+1}) = a_{k+1,k} = q_{v_{k+1}-v_k} = q_{v_k-v_{k+1}}
= (M(\Burau{\sigma_k}v,\Burau{\sigma_k}v))_{k,k+1}
\]
and similarly $\hom{\sigma_k}(a_{k+1,k}) =
(M(\Burau{\sigma_k}v,\Burau{\sigma_k}v))_{k+1,k}$; and finally, if
$i,j \neq k,k+1$, then $\hom{\sigma_k}(a_{ij}) = q_{v_i-v_j} =
(M(\Burau{\sigma_k}v,\Burau{\sigma_k}v))_{ij}$. This shows that
$\hom{\sigma_k}(A) = M(\Burau{\sigma_k}v,\Burau{\sigma_k}v)$.

By replacing $v$ by $\Burau{\sigma_k^{-1}}v$, we can deduce the
desired identity for $B=\sigma_k^{-1}$ from the identity for
$B=\sigma_k$. It thus suffices to show that if the identity holds
for $B=B_1$ and $B=B_2$, then it holds for $B=B_1B_2$. Now
\begin{align*}
\hom{B_1B_2}(M(v,v)) &=
\hom{B_1}(\hom{B_2}(M(v,v))) \\
&=
\hom{B_2}(M(\Burau{\hat{B_1}}v,\Burau{\hat{B_1}}v)) \\
&=
M(\Burau{\hat{B_2}}\Burau{\hat{B_1}}v,\Burau{\hat{B_2}}\Burau{\hat{B_1}}v) \\
&= M(\Burau{\widehat{B_1B_2}}v,\Burau{\widehat{B_1B_2}}v),
\end{align*}
and the lemma follows.
\end{proof}

\begin{lemma}
Let $B\in B_n$ and $v\in\Z^n$. If $A=M(v,v)$, then
$\Phil{B}(A) \cdot A = M(\Burau{\hat{B}}v,v)$
and $A \cdot \Phir{B}(A) = M(v,\Burau{\hat{B}}v)$.
\label{qlemma2}
\end{lemma}

\begin{proof}
We will prove the identity involving $\Phil{}$; the proof for
$\Phir{}$ is completely analogous. First note that the identity is
obvious when $B$ is trivial. Now suppose that the identity holds for
$B$; we show that it holds for $B\sigma_k$. Indeed, by
Proposition~\ref{chainrule}, we have
\[
\Phil{B\sigma_k}(A) \cdot A = \Phil{\sigma_k}(\hom{B}(A)) \cdot \Phil{B}(A) \cdot A
= \Phil{\sigma_k}(\hom{B}(A)) \cdot M(\Burau{\hat{B}}v,v).
\]
By Lemma~\ref{sigmak}, $(\Phil{B\sigma_k}(A) \cdot A)_{ij} =
(\Burau{\hat{B}}v)_i-v_j$ except when $i=k$ or $i=k+1$; in these
cases, Lemmas~\ref{sigmak}, \ref{funceq}, and \ref{qlemma1} yield
\begin{align*}
(\Phil{B\sigma_k}(A) \cdot A)_{kj} &= -\hom{B}(a_{k+1,k})
(M(\Burau{\hat{B}}v,v))_{kj}
- (M(\Burau{\hat{B}}v,v))_{k+1,j} \\
&= -(q_{(\Burau{\hat{B}}v)_{k+1} - (\Burau{\hat{B}}v)_k})
(q_{(\Burau{\hat{B}}v)_k-v_j})
- q_{(\Burau{\hat{B}}v)_{k+1}-v_j} \\
&= q_{(2(\Burau{\hat{B}}v)_k - (\Burau{\hat{B}}v)_{k+1} - v_j)}, \\
(\Phil{B\sigma_k}(A) \cdot A)_{k+1,j} &=
(M(\Burau{\hat{B}}v,v))_{kj} = q_{(\Burau{\hat{B}}v)_k - v_j}.
\end{align*}
In all cases, we have $(\Phil{B\sigma_k}(A) \cdot A)_{ij} =
q_{(\Burau{\sigma_k \hat{B}}v)_i - v_j}$, and so the identity holds
for $B\sigma_k$. A similar computation shows that if the identity
holds for $B$, then it holds for $B\sigma_k^{-1}$. The identity for
all $B$ follows by induction.
\end{proof}

\begin{proof}[Proof of Theorem \ref{detthm}]
Since $\red{\hat{B}}$ is a presentation matrix for
$H_1(\Sigma_2(K))=H_1(\Sigma_2(\overline{K}))$, and $n(K)$ is by
definition an invariant factor of $H_1(\Sigma_2(K))$, there exists a
vector $v=(v_2,\dots,v_n)\in\Z^{n-1}$ with
$\gcd(v_2,\dots,v_n,n(K))=1$, such that the vector
$(\red{\hat{B}}-1)v$ is divisible by $n(K)$. This means that
$v'=(0,v_2,\dots,v_n)\in\Z^n$ satisfies the following property: if
we write $w=(\Burau{\hat{B}}-1)v'$, then $w_i \equiv w_j \pmod
{n(K)}$ for all $1\leq i,j\leq n$. Since the Burau representation
sends $(1,\dots,1)$ to itself, we may define
$v''=(-w_1,v_2-w_1,\dots,v_n-w_1)$ and conclude that the vector
$(\Burau{\hat{B}}-1)v''$ is divisible by $n(K)$.

Now consider the ring homomorphism $\rho\co \A_n \to \Z[x]$ defined
by sending the matrix $A$ to $M(v'',v'')$.  By Lemma~\ref{qlemma2},
$\rho$ sends $(1-\Phil{B}(A)) \cdot A$ to
$M(v'',v'')-M(\Burau{\hat{B}}v'',v'')$. Since $\Burau{\hat{B}}v''
\equiv v'' \pmod {n(K)}$, Lemma~\ref{divisible} shows that for all
$i,j$,
\[
p_{(n(K)+1)/2} \,\left|\, q_{v''_i-v''_j} -
q_{(\Burau{\hat{B}}v'')_i-v''_j} =
\left(M(v'',v'')-M(\Burau{\hat{B}}v'',v'')\right)_{ij} \right. .
\]
Hence $p_{(n(K)+1)/2}$ divides all entries of $\rho((1-\Phil{B}(A))
\cdot A)$. Similarly, $p_{(n(K)+1)/2}$ divides all entries of
$\rho(A \cdot (1-\Phir{B}(A)))$. It follows that $\rho$ descends to
a map $\tilde{\rho}\co  HC_0(K) \to \Z[x]/(p_{(n(K)+1)/2}(x))$,
which clearly factors through $HC_0^{\ab}(K)$ since $q_m = q_{-m}$.

It remains to show that $\tilde{\rho}$ is surjective, or
equivalently, that $x$ is in the image of $\tilde{\rho}$. Define the
set $S = \{j \in\Z\,|\,q_j(x) \in \im \rho\}$. We first claim that
$S = g\Z$ for some $g\in\Z$.

By Lemma~\ref{funceq}, if $m_1\in S$, then $m_1+m_2\in S$ if and
only if $m_1-m_2\in S$. Hence if $a_1,a_2\in S$, then the arithmetic
progression $\{a_1+\ell(a_2-a_1)\}_{\ell\in\Z}$ is contained in $S$.
Now suppose that $(m_1,m_2)=(a_1,a_2)$ minimizes $|m_1-m_2|$, where
$(m_1,m_2)$ ranges over all pairs of distinct elements in $S$. Then
$S=\{a_1+\ell(a_2-a_1)\}_{\ell\in\Z}$, since the existence of an
element of $S$ not in the arithmetic progression would violate
minimality. Furthermore, $0\in S$ since $q_0(x) = -2 \in \im \rho$.
Hence $S$ is an ideal in $\Z$, which proves the claim.

Now for $2\leq i\leq n$, we have $\rho(a_{i1}) = q_{v''_i-v''_1} =
q_{v_i}$, and so $v_i \in S$. By the claim, we conclude that all
multiples of $\gcd(v_2,\dots,v_n)$ are in $S$. Since
$\gcd(v_2,\dots,v_n)$ and $n(K)$ are relatively prime by
construction, there exists $j\in S$ with $j\equiv 1 \pmod {n(K)}$.
By Lemma~\ref{divisible}, $p_{(n(K)+1)/2}$ divides $q_j(x)-q_1(x)$;
but $q_j(x) \in \im \rho$ and $q_1(x) = -x$, and hence $x \in \im
\tilde{\rho}$, as desired.
\end{proof}

The above proof is a bit tedious and hard to digest. In \cite{II},
we will use the cord ring to see the surjection from
Theorem~\ref{detthm} rather more directly, and justify the
appearance of the polynomials $p_m$.

%*********************************************************************
\subsection{Linearized knot contact homology}
\label{ssec:lin}

Although it is generally difficult to find a nice form for the full
knot contact homology $HC_*(K)$, there is a canonical linearization
of $HC_*(K)$ which is easy to compute, and which encodes the group
$H_1(\Sigma_2(K))$ and hence the determinant of $K$. We will discuss
this linearization in this section.

Let $(\A,\d)$ be the knot DGA for a knot $K$. It follows easily from
Proposition~\ref{prop:augment} that the algebra map $\epsilon\co  \A
\to \Z$, sending each generator $a_{ij}$ to $-2$, and all other
generators to $0$, gives an augmentation for $(A,\d)$; that is,
$\epsilon \circ \d = 0$. As in \cite{Che}, we can construct
linearized homology groups via this augmentation.

Let $\varphi_\epsilon\co  \A \to \A$ be the algebra isomorphism
which sends $a_{ij}$ to $a_{ij}-2$ for all $i,j$, and acts as the
identity on the other generators $b_{ij},c_{ij},d_{ij},e_i$. Let
$\M$ be the subalgebra of $\A$ generated by all words of length at
least $1$. Then $\d_\epsilon := \varphi_\epsilon \circ \d \circ
\varphi_\epsilon^{-1}$ maps $\M$ into itself. For any $k\geq 2$,
$\d_\epsilon$ descends to a differential $\d_\epsilon^{\lin}$ on the
finite-rank graded $\Z$-module $\M/\M^k$, and the homology of the
resulting complex gives a linearized version of $HC_*(K)$.

We now restrict our attention to the simplest nontrivial case,
$k=2$, although higher values of $k$ also yield knot invariants.
Note that $\M/\M^2$ is a free $\Z$-module generated in degree $0$ by
$a_{ij}$, in degree $1$ by $b_{ij}$ and $c_{ij}$, and in degree $2$
by $d_{ij}$ and $e_i$.

\begin{definition}
The \textit{linearized contact homology} of $K$, written
$HC_*^{\lin}(K)$, is the graded homology of $\d_{\epsilon}$ on
$\M/\M^2$.
\end{definition}

Using the conjugation symmetry of Section~\ref{ssec:conjugation}, we
can also define an abelian version of linearized contact homology.
In this case, $(\A,\d)$ is the abelian knot DGA of $K$ with
corresponding maximal subalgebra $\M$, and $\d_{\epsilon}$ descends
to a differential on $\M/\M^2$, whose homology we denote by
$HC_*^{\ab,\lin}(K)$.

\begin{proposition}
The groups $HC_*^{\lin}(K)$ and $HC_*^{\ab,\lin}(K)$ are invariants
of the knot $K$.
\end{proposition}

\begin{proof}
All of the stable tame isomorphisms used in the proof of
Theorem~\ref{knotthm} commute with $\epsilon$, because $\phi_B \circ
\epsilon = \epsilon \circ \phi_B$ for any braid $B$. The proposition
follows.
\end{proof}

We now establish some results about $HC_*^{\lin}(K)$. Note that this homology exists
only in dimensions $* = 0,1,2$.

\begin{proposition}
For $K$ a knot, $HC_2^{\lin}(K) \cong \Z$ and
$HC_1^{\lin}(K) \cong (H_1(\Sigma_2) \oplus \Z) \otimes (H_1(\Sigma_2) \oplus \Z)
\oplus \Z^m$ for some $m \geq 0$.
\label{prop:HClin1}
\end{proposition}

Proposition~\ref{prop:HClin1} will follow almost immediately from a lemma
which we now formulate. Let $d_1,d_2,\ldots,d_k$ be the invariant factors of
$H_1(\Sigma_2)$, with $d_1|d_2|\cdots|d_k$. Write $Y = \Bur{\hat{B}}$,
and express $1-Y$ in Smith normal form; more precisely, there are matrices
$X_1,X_2 \in GL(n,\Z)$ with $X_1\cdot (1-Y) X_2 = \Delta$, where $\Delta$
is a diagonal matrix. Since $1-Y$ is a presentation matrix for
$H_1(\Sigma_2) \oplus \Z$, we may assume that the diagonal entries of
$\Delta$ are precisely $d_1,\ldots,d_n$, where $d_{k+1} = \cdots = d_{n-1} = 1$
and $d_n = 0$.

Now write the linearized chain complex of the alternate knot DGA of
the braid $B$ as $C_2 \stackrel{M_2}{\rightarrow} C_1
\stackrel{M_1}{\rightarrow} C_0$, where $C_0$, $C_1$, $C_2$ are
generated by $\{a_{ij}\}$, $\{b_{ij},c_{ij},d_{ij}\}$,
$\{e_{ij},f_{ij}\}$, respectively. The map $M_2$ is given by the
action of $\d_\epsilon^{\lin}$ on $C_2$: $\d_\epsilon^{\lin}E = B -
D - C \cdot (\Bur{\hat{B}})^{\transpose}$, $\d_\epsilon^{\lin}F = B
- C - \Bur{\hat{B}} \cdot D$.

\begin{lemma}
After changing bases for $C_1$ and $C_2$, we can write $M_2$ as a diagonal
$(3n^2-n) \times (2n^2)$ matrix, with $n^2$ diagonal entries given by
$\gcd(d_i,d_j) = d_{\min(i,j)}$, $1 \leq i,j \leq n$, and the remaining
$n^2$ diagonal entries given by $1$.
\label{lem:M2}
\end{lemma}

\begin{proof}
We wish to calculate Smith normal form for $M_2$. Replace $D$ by
$D+B-C\cdot Y^{\transpose}$ and then $F$ by $F+Y\cdot E$ (this
simply reparametrizes $d_{ij}$ and $f_{ij}$); then
$\d_\epsilon^{\lin}E = -D$ and $\d_\epsilon^{\lin}F = B - C - Y
\cdot (B-C\cdot Y^{\transpose})$. The action of $\d_\epsilon^{\lin}$
on $E$ gives $n^2$ ``$1$'' diagonal entries in $M_2$. To complete
the Smith normal form for $M_2$, we need to examine the matrix $Z =
B - C - Y \cdot (B-C\cdot Y^{\transpose})$.

Introduce dummy variables $b_{11},\ldots,b_{nn}$ (this adds $n$ rows
of zeros to $M_2$), and replace $b_{ij}$ by $b_{ij}-b_{jj}$ for all
$i\neq j$. If we write $\tilde{B} = (b_{ij})$, which is identical to
$B$ away from the diagonal, then $Z$ becomes $\tilde{B} - C - Y
\cdot (\tilde{B}-C\cdot Y^{\transpose})$. Next replace $Z$ by $X_1
\cdot Z \cdot Z_1^{\transpose}$; this has the effect of changing
basis in $C_2$. Finally, change basis in $C_1$ by successively
replacing $\tilde{B}$ by $X_2 \cdot \tilde{B} \cdot
(X_1^{\transpose})^{-1} +C$, and then $C$ by $Y^{-1}\cdot X_1^{-1}
\cdot C \cdot X_2^{\transpose}$. Then $Z$ becomes
\[
\Delta \cdot \tilde{B} - C \cdot \Delta = (d_i b_{ij} - d_j c_{ij}).
\]
The lemma follows immediately.
\end{proof}

\begin{proof}[Proof of Proposition~\ref{prop:HClin1}]
We have $HC_2^{\lin}(K) = \ker M_2$ and $HC_1^{\lin}(K) = (\ker M_1)/(\im M_2)$.
Write $M_2$ as in Lemma~\ref{lem:M2}. Since $d_1,\ldots,d_{n-1}>0$,
$M_2$ has only one $0$ diagonal entry, and so $\ker M_2 \cong \Z$.

For $HC_1^{\lin}(K)$, assume that we have chosen bases for $C_1,C_2$
so that $M_2$ is diagonal with diagonal entries
$\delta_1,\ldots,\delta_{2n^2}$, with $\delta_{2n^2}=0$ and
$\delta_i \neq 0$ for $i<2n^2$. Let the corresponding basis for
$C_1$ be $v_1,\ldots,v_{3n^2-n}$. Then $v_1,\ldots,v_{2n^2-1} \in
\ker M_1$, and indeed we may assume that
$\{v_1,\ldots,v_{2n^2-1+m}\}$ forms a basis for $\ker M_1$ for some
$m \geq 0$. It follows that $(\ker M_1)/(\im M_2) \cong \Z^m \oplus
\oplus_i \Z_{\delta_i}$. From Lemma~\ref{lem:M2}, we conclude that
$(\ker M_1)/(\im M_2) \cong (H_1(\Sigma_2) \oplus \Z) \otimes
(H_1(\Sigma_2) \oplus \Z) \oplus \Z^m$, as desired.
\end{proof}

We conjecture that $m$ is always $0$ in
Proposition~\ref{prop:HClin1}. Even if this is not the case,
Proposition~\ref{prop:HClin1} still implies that we can deduce
$H_1(\Sigma_2)$ from $HC_1^{\lin}(K)$; this is an easy consequence
of the fundamental theorem of finitely generated abelian groups.

\begin{corollary}
The group $H_1(\Sigma_2)$, and hence the determinant
$|\Delta_K(-1)|$, is determined by the equivalence class of the knot
DGA of $K$.
\label{cor:determined}
\end{corollary}

We now turn our attention to $HC_0^{\lin}(K)$, and derive a result
which will be used in the sequel \cite{II}.

\begin{proposition}
\label{prop:linsurject}
There are surjections
\[
HC_0^{\lin}(K) \twoheadrightarrow HC_0^{\lin,\ab}(K)
\twoheadrightarrow \Sym^2(H_1(\Sigma_2)).
\]
\end{proposition}

The proof of Proposition~\ref{prop:linsurject} is essentially a
linearized version of the (more involved) proof of
Theorem~\ref{detthm}, and we provide just an outline here. The
matrix $\Bur{\hat{B}}$ acts on $\Z^n$, which has generators
$e_1,\ldots,e_n$, and it is an easy exercise to check that
$H_1(\Sigma_2)$ is generated by the differences $\{e_i-e_j\}$,
modulo the relations $\Bur{\hat{B}}e_i - e_i = 0$ for all $i$. (Note
that $\Bur{\hat{B}}e_i-e_i$ is a linear combination of differences
$e_{j_1}-e_{j_2}$.)

\begin{lemma}
\label{lem:linsurject}
In the ring $\Z[e_1,\ldots,e_n]$, if $A = (-2+(e_i-e_j)^2)$, then
\[
\Phil{B}(A) \cdot A = \left(-2+(\Bur{\hat{B}}e_i-e_j)^2 + O(e^4)\right)
\]
and
\[
\Phir{B}(A) \cdot A = \left(-2+(e_i-\Bur{\hat{B}}e_j)^2 + O(e^4)\right).
\]
\end{lemma}

\begin{proof}
Induction on the braid word $B$, as in the proof of Lemma~\ref{qlemma2}.
\end{proof}

\begin{proof}[Proof of Proposition~\ref{prop:linsurject}]
The surjection from $HC_0^{\lin}(K)$ to $\Sym^2(H_1(\Sigma_2))$
which factors through $HC_0^{\lin,\ab}(K)$, sends $a_{ij}$ to
$(e_i-e_j)^2$. This is well-defined by Lemma~\ref{lem:linsurject}
and the fact that $\Bur{\hat{B}}e_i-e_i=0$ in $H_1(\Sigma_2)$ for
all $i$; it is a surjection by polarization, since $2$ is invertible
in $\Sym^2(H_1(\Sigma_2))$.
\end{proof}

A more natural way to see the second surjection in
Proposition~\ref{prop:linsurject}, using cords, is given in
Section~5.2 of \cite{II}. By all computational indications, it seems
likely that this surjection is always an isomorphism, and hence that
$HC_0^{\lin,\ab}(K)$ is determined by $H_1(\Sigma_2)$ for all knots
$K$. However, $HC_0^{\lin}(K)$ is not determined by $H_1(\Sigma_2)$.
For instance, the knots $8_{18}$ and $9_{37}$ both satisfy
$H_1(\Sigma_2) \cong \Z_3 \oplus \Z_{15}$, but $HC_0^{\lin}(8_{18})
\cong \Z_3^3 \oplus \Z_{15}$, while $HC_0^{\lin}(9_{37}) \cong
\Z_3^2 \oplus \Z_{15}$. It is possible that $HC_0^{\lin}(K)$ is
always isomorphic to either $\Sym^2(H_1(\Sigma_2))$ or
$H_1(\Sigma_2) \otimes H_1(\Sigma_2)$; this has been verified for
all prime knots with eleven or fewer crossings.

%*********************************************************************
%*********************************************************************
\section{Computations}
\label{sec:computations}

The combinatorial nature of our invariants allows them to be
calculated readily by computer. \textit{Mathematica} source to
compute the invariants is available at the author's web site. The
full program requires the installation of the noncommutative algebra
package NCAlgebra/NCGB \cite{NCAlg}, though there is a commutative
version which does not use NCAlgebra. Using the program, we can
calculate the braid homomorphism $\phi$, braid and knot DGAs, $HC_0$
(via Gr\"obner bases), augmentation numbers, and linearized knot
contact homology. Many of the calculations in this section were
performed with the help of the program.

%*********************************************************************
\subsection{Braid conjugacy classes}
\label{ssec:braidcomp}

It can be difficult to tell if two braids have equivalent braid
DGAs. Currently, the most effective computational tool in this
regard is the augmentation number invariant described in
Section~\ref{ssec:augment}, which can be calculated by computer and
is reasonably effective in distinguishing between braid conjugacy
classes.

\begin{table}[ht!]\small
\[
\begin{array}{|c||c|c|c|c|c|}
\hline \heighten B & w & \Aug(B,2) & \Augab(B,2) & \Augab(B,5) & \Augab(B,7) \\
\hline\hline
\heighten\sigma_1^5\sigma_2 & 6 & 5 & 1 & 4 & 6\\
\hline
\heighten\sigma_1^3\sigma_2^3 & 6 & 9 & 5 & 12 & 18\\
\hline
\heighten\sigma_1^3\sigma_2 & 4 & 4 & 2 & 5 & 7\\
\hline
\heighten\sigma_1^5\sigma_2^{-1} & 4 & 4 & 2 & 5 & 6\\
\hline
\heighten\sigma_1^3\sigma_2\sigma_1^{-1}\sigma_2 & 4 & 4 & 2 & 4 & 7\\
\hline
\heighten\sigma_1\sigma_2 & 2 & 4 & 2 & 5 & 7\\
\hline
\heighten\sigma_1^3\sigma_2^{-1} & 2 & 4 & 2 & 4 & 7\\
\hline
\heighten\sigma_1^3\sigma_2^{-1}\sigma_1\sigma_2^{-1} & 2 & 2 & 2 & 9 & 4\\
\hline
\heighten\sigma_1\sigma_2^{-1} & 0 & 5 & 1 & 4 & 6\\
\hline
\heighten\sigma_1\sigma_2^{-1}\sigma_1\sigma_2^{-1} & 0 & 9 & 5 & 12 & 18\\
\hline
\heighten\sigma_1^3\sigma_2^{-3} & 0 & 7 & 5 & 11 & 13\\
\hline
\sigma_1^2\sigma_2^{-2}\sigma_1\sigma_2^{-1} & 0 & 1 & 1 & 4 & 2\\
\hline
\heighten\sigma_1^{-2}\sigma_2^2\sigma_1^{-1}\sigma_2 & 0 & 1 & 1 & 4 & 2\\
\hline
\end{array}
\]
\caption{
The conjugacy classes in $B_3$ with nonnegative writhe which can be
represented by a connected braid of word length at most $6$,
identified by such a representative $B$. The writhe $w$ is given,
along with various augmentation numbers.
}
\label{braidtable}
\end{table}

In Table~\ref{braidtable}, we list all conjugacy classes in $B_3$ of
nonnegative writhe which can be represented by a connected braid of
word length at most $6$, along with their writhe and selected
augmentation numbers. (Recall that the writhe $w$ is the
homomorphism $B_n \to \Z$ for which $w(\sigma_k) = 1$ for all $k$.)
%this obviously descends to an invariant of braid conjugacy class.
The augmentation numbers, combined with the writhe, distinguish all
of the given conjugacy classes, with the exception of the mirrors
$\sigma_1^2\sigma_2^{-2}\sigma_1\sigma_2^{-1}$ and
$\sigma_1^{-2}\sigma_2^2\sigma_1^{-1}\sigma_2$, which by
Proposition~\ref{braidmirror} cannot be distinguished using our
invariants.

Without the use of the writhe, our invariants become less effective.
For instance, it can be shown that the braid DGAs for
$\sigma_1\sigma_2$ and $\sigma_1^3\sigma_2$ are both stable tame
isomorphic to the DGA on two generators of degree $0$ and two
generators of degree $1$, with $\d=0$.

%*********************************************************************
\subsection{The unknot and unlink, and other links}
\label{ssec:unknot}

The unknot is the closure of the trivial braid in $B_1$, whose knot
DGA is generated by $b_{11},c_{11}$ in degree $1$ and $d_{11},e_1$
in degree $2$, with differential given by $\d b_{11} = \d c_{11} =
\d d_{11} = 0$, $\d e_1 = b_{11}+c_{11}$. (The $1\times 1$ matrices
$\Phil{}$ and $\Phir{}$ for the trivial braid in $B_1$ are simply
$(1)$.) This is stable tame isomorphic to the DGA with generators
$b_{11}$ in degree $1$ and $d_{11}$ in degree $2$, and with $\d =
0$. In particular, we see that $HC_0(\textrm{unknot}) \cong \Z$.

More generally, let $L_n$ denote the $n$-component unlink, which is
the closure of the trivial braid $B^0_n$ in $B_n$. The matrices
$\Phil{B^0_n}$ and $\Phir{B^0_n}$ are both the identity, and so the
knot DGA of $B^0_n$ has differential $\d a_{ij} = \d b_{ij} = \d
c_{ij} = \d d_{ij} =0$, $\d e_i = b_{ii}+c_{ii}$, and $HC_0(L_n)
\cong \Z\langle \{a_{ij}\}\rangle$.

By contrast, the Hopf link, which is the closure of $\sigma_1^2 \in
B_2$, has knot DGA satisfying
$$\d B = \left( \begin{matrix} -4+a_{12}a_{21} & 4a_{12}-a_{12}a_{21}a_{12} \\
0 & -4+a_{21}a_{12} \end{matrix} \right)$$
$$\d C = \left( \begin{matrix} -4+a_{12}a_{21} & 0 \\
4a_{21}-a_{21}a_{12}a_{21} & -4+a_{21}a_{12} \end{matrix} \right).
\leqno{\rm and}$$
It follows that $HC_0(\textrm{Hopf}) \cong \Z\langle a_{12},a_{21}\rangle / \langle -4+a_{12}a_{21},
-4+a_{21}a_{12} \rangle$, which is not isomorphic
to $HC_0(L_2)$.

Finally, consider the split link $3_1\sqcup 0_1$ given by the
closure of $\sigma_1^3\in B_3$, which is the unlinked union of an
unknot and a trefoil. We can compute that
\begin{align*}
HC_0(3_1\sqcup 0_1) &\cong \Z\langle
a_{12},a_{13},a_{31},a_{23},a_{32}\rangle /  \langle
-2+a_{12}+a_{12}^2,
\\
& \qquad a_{13}-a_{23}-a_{12}a_{13}+a_{12}a_{23},
a_{31}-a_{32}-a_{31}a_{12}+a_{32}a_{12} \rangle,
\end{align*}
which is not isomorphic to $HC_0(L_2)$. Hence our invariants, in
contrast to the Alexander polynomial, can distinguish between split
links.

%*********************************************************************
\subsection{Torus knots $T(2,n)$}
\label{ssec:T2}

In this section, we generalize our computation in
Section~\ref{ssec:HCdef} of $HC_0$ for the trefoil, by computing
$HC_0$ for all torus knots of the form $T(2,n)$.

First we introduce two sequences related to the polynomials $q_m(x)$
discussed in Section~\ref{ssec:det}. In $\A_2$, define
$\{q_m^{(1)}\}$, $\{q_m^{(2)}\}$ by $q_0^{(1)}=q_0^{(2)}=-2$,
$q_1^{(1)}=a_{12}$, $q_1^{(2)}=a_{21}$, and $q_{m+1}^{(1)} =
-a_{12}q_m^{(2)}-q_{m-1}^{(1)}$, $q_{m+1}^{(2)} =
-a_{21}q_m^{(1)}-q_{m-1}^{(2)}$; by definition, we have
$q_m^{(1)}|_{a_{12}=a_{21}=-x} = q_m^{(2)}|_{a_{12}=a_{21}=-x} =
q_m(x)$.

\begin{lemma}
For $k \geq 1$, we have \label{torusknot}
\[
\Phil{\sigma_1^{2k-1}}(A) \cdot A = \left(
\begin{matrix}
q_{2k-1}^{(2)} & q_{2k}^{(2)} \\
q_{2k-2}^{(1)} & q_{2k-1}^{(1)}
\end{matrix} \right)
\quad \textrm{and} \quad
A \cdot \Phir{\sigma_1^{2k-1}}(A) = \left(
\begin{matrix}
q_{2k-1}^{(1)} & q_{2k-2}^{(1)} \\
q_{2k}^{(2)} & q_{2k-1}^{(2)}
\end{matrix} \right) .
\]
\end{lemma}

\begin{proof}
We prove the first identity; the second follows similarly or by
conjugation symmetry. For $k=1$, the identity can be verified
directly. For $k \geq 2$, by Proposition~\ref{chainrule}, we have
\[
\Phil{\sigma_1^{2k-1}}(A) \cdot A
= \Phil{\sigma_1^2}(\hom{\sigma_1^{2k-3}}(A)) \cdot \Phil{\sigma_1^{2k-3}}(A) \cdot A.
\]
Now $\hom{\sigma_1^2}(A) = A$, and so
$\hom{\sigma_1^{2m-1}}(A) = \hom{\sigma_1}(A) =
\left( \begin{smallmatrix} -2 & a_{21} \\ a_{12} & -2 \end{smallmatrix} \right)$;
then we compute that $\Phil{\sigma_1^2}(\hom{\sigma_1^{2m-1}}(A)) =
\left( \begin{smallmatrix} -1+a_{21}a_{12} & a_{21} \\ -a_{12} & -1 \end{smallmatrix} \right)$.
The desired identity can now be checked by induction on $m$.
\end{proof}

We can now compute $HC_0$ for the $(2,2k-1)$ torus knot.

\begin{proposition}
For $k \geq 1$, $HC_0(T(2,2k-1)) \cong \Z[x]/(p_k(x))$, where
$p_k$ is the polynomial defined in Section~\ref{ssec:det}.
\label{T2nprop}
\end{proposition}

\begin{proof}
By Lemma~\ref{torusknot}, we have
$((1-\Phil{\sigma_1^{2k-1}}(A))\cdot A)_{12} = a_{12} -
q_{2k}^{(2)}$ and $(A\cdot (1-\Phir{\sigma_1^{2k-1}}(A)))_{21} =
a_{21} - q_{2k}^{(2)}$. Hence in $HC_0(T(2,2k-1))$, we have
$a_{12}=a_{21}$, and we can set $x:=-a_{12}=-a_{21}$. The relations
in $HC_0(T(2,2k-1))$ then become the entries of the matrix
$$\left( \begin{matrix}
-2-q_{2k-1}(x) & -x-q_{2k}(x) \\
-x-q_{2k-2}(x) & -2-q_{2k-1}(x)
\end{matrix}
\right).$$
Now define $r_m(x) = q_m(x)-q_{2k-1-m}(x)$ for all $m$, so that
$r_m(x) = xr_{m+1}(x) - r_{m+2}(x)$. We have
\begin{align*}
HC_0(T(2,2k-1)) &\cong
\Z[x]/(-x-q_{2k-2}(x),-x-q_{2k}(x),-2-q_{2k-1}(x)) \\
&= \Z[x]/(r_0(x),r_1(x)),
\end{align*}
where the second equality holds because $-x-q_{2k}(x)= r_{-1}(x) =
xr_0(x)-r_1(x)$. Now
\[
\gcd(r_0(x),r_1(x)) = \gcd(r_1(x),r_2(x)) = \cdots =
\gcd(r_{k-1}(x),r_k(x)) = r_{k-1}(x);
\]
since $r_{k-1}(x)=p_k(x)$ by induction on $k$, the proposition
follows.
\end{proof}

Proposition~\ref{T2nprop} shows that the surjection in
Theorem~\ref{detthm} is actually an isomorphism for $(2,n)$ torus
knots. This is established in more generality, for all two-bridge
knots, in \cite{II}. However, we will see in the following sections
that it is not true for all knots.

Before moving on, we make an observation regarding the trefoil in
the context of Legendrian knot theory. We have seen that $HC_0$ for
the trefoil is $\Z[x]/(x^2+x-2)$. This has two maps to $\Z_2$ and
thus the knot DGA for the trefoil has two augmentations over $\Z_2$.
If we use the form for the trefoil knot DGA given at the end of
Section~\ref{ssec:invariants}, then one augmentation, $\epsilon_0$,
sends $a_{12}$ and $a_{21}$ to $0$, while the other, $\epsilon_1$,
sends them to $1$. The linearized differential for $\epsilon_0$ is
simply given by the linear terms in the differential mod $2$, and
one can calculate that this has Poincar\'e polynomial
$\lambda^2+\lambda$ (ie, homology of dimension $1$ in degrees $2$
and $1$). On the other hand, one can check that the augmentation
$\epsilon_1$ yields Poincar\'e polynomial $4\lambda^2+4\lambda$.
This is an example where different augmentations of a DGA measuring
relative contact homology give different Poincar\'e polynomials over
$\Z_2$; this phenomenon has also now been observed for some
Legendrian knots in $\R^3$ \cite{MS}.

%*********************************************************************
\subsection{Torus knots $T(3,n)$}
\label{ssec:T3n}

The invariant $HC_0$ is not as simple for general torus knots as it is
for $T(2,n)$. For instance, if we write the torus knot $T(3,4)$ as the closure
of the braid $(\sigma_1\sigma_2)^4 \in B_3$, we can calculate
\begin{align*}
HC_0(T(3,4)) &\cong \Z\langle a_1,a_2\rangle / \langle
2a_1-2a_2-a_1^2+a_1a_2+a_1^3-a_1^2a_2, \\
& \qquad -a_1a_2+a_2a_1+a_1^2a_2-a_2a_1^2,
-3a_1+3a_2+a_1a_2a_1-a_2a_1a_2,
\\
& \qquad -a_1+a_2+a_1^2-a_2^2, 2+a_2-a_1a_2-a_2a_1-a_1^3 \rangle .
\end{align*}
Furthermore, if we use lexicographic order with $a_1<a_2$, then the
above presentation gives a Gr\"obner basis for the ideal. Thus as an
additive abelian group, $HC_0(T(3,4)) \cong \Z^7$, generated by
$1,a_1,a_2,a_1^2,a_1a_2,a_2a_1,a_1a_2a_1$. On the other hand, the
abelianization of $HC_0(T(3,4))$ is easily shown to be isomorphic to
$\Z[x]/(x^3-2x^2-x+2)$, and so $HC_0(T(3,4))$ is noncommutative.
Note that $HC_0(T(3,4))$ does map surjectively to $\Z[x]/(x^2-x-2)$,
as per Theorem~\ref{detthm}.

In general, it seems that $HC_0(T(m,n))$ is too complicated to be
expressed in closed form. However, for $T(3,n)$, there is a
relatively simple expression for $HC_0^{\ab}$, if not for $HC_0$
(cf\ $HC_0(T(3,4))$ above).

It is easy to see that $HC_0^{\ab}(T(3,n))$ can be written in the
form $\Z[x]/(p(x))$ for some polynomial $p$. The knot $T(3,n)$ is
the closure of $(\sigma_1\sigma_2)^n \in B_3$, and
$\phi_{(\sigma_1\sigma_2)^n}$ cyclically permutes
$a_{12},a_{23},a_{31}$ when $3$ does not divide $n$. Hence the
generators $a_{12},a_{23},a_{13}$ of $HC_0^{\ab}(T(3,n))$ are equal
by Corollary~\ref{quotientcor}. (The same argument shows that
$HC_0(T(3,n))$ is generated by two elements.)

Explicitly calculating $HC_0^{\ab}(T(3,n))$ is more involved; a
straightforward but tedious induction, along the lines of the proof
of Proposition~\ref{T2nprop}, yields the following result.

\begin{proposition}
If $\gcd(3,n)=1$, then
\[
HC_0^{\ab}(T(3,n)) \cong \begin{cases}
\Z[x]/(r_k(x)), & n=2k-1 \\
\Z[x]/(s_k(x)), & n=2k-2,
\end{cases}
\]
where $r_0(x)=2-x$, $r_1(x)=x-2$, $r_{m+1}(x) = (x-1)r_m(x)-r_{m-1}(x)$, and
$s_1(x)=0$, $s_2(x)=(x-2)(x+1)$, $s_{m+1}(x) = (x-1)s_m(x)-s_{m-1}(x)$.
\label{T3nprop}
\end{proposition}

\noindent Note that $r_k(x) = \frac{x-2}{x-3} \, p_k(x-1)$.

%*********************************************************************
\subsection{The connected sums $3_1 \# 3_1$ and $3_1 \# \overline{3_1}$}
\label{ssec:connectsum}

When we view the connected sum $3_1 \# 3_1$ as the closure of
$\sigma_1^3\sigma_2^3 \in B_3$, we can calculate that $\Augab(3_1\#
3_1,2)=7$, and hence $HC_0(3_1\# 3_1)$ is not of the form
$\Z[x]/(p(x))$.  To garner more information, we use the standard
Gr\"obner basis algorithm to find that
\begin{align*}
HC_0(3_1 \# 3_1) &\cong \Z\langle a_1,a_2,a_3,a_4 \rangle / \langle
(a_1+2)(a_1-1), (a_2+2)(a_2-1), \\
& \qquad (a_1-1)a_3, a_4(a_1-1), (a_2-1)a_4, a_3(a_2-1) \rangle .
\end{align*}
The given presentation for the above ideal is a noncommutative
Gr\"obner basis for the ideal with respect to any lexicographic
order, and it follows that $HC_0(3_1 \# 3_1)$ is noncommutative (for
instance, $a_1a_3 \neq a_3a_1$) and not finitely generated as a
$\Z$-module. It is also true that $HC_0^{\ab}$ has infinite rank in
this case; in fact,
\begin{align*}
HC_0^{\ab}(3_1 \# 3_1) &\cong \Z[a_1,a_2,a_3] / ((a_1+2)(a_1-1),
(a_2+2)(a_2-1), \\
& \qquad (a_1-1)(a_3-a_2), (a_2-1)(a_3-a_1)).
\end{align*}

As another application of this example, it turns out that precisely
the same computation holds for $3_1\#\overline{3_1}$; hence it is
not true that $HC_0$ or $HC_0^{\ab}$ determines the isotopy class of
a knot, even up to mirrors.

%*********************************************************************
\subsection{Knot contact homology and other invariants}
\label{ssec:HOMFLY}

From the computations of the previous sections, the $HC_0$ invariant
for knots depends on more than the determinant of the knot. It is
also not determined by the Alexander polynomial; for instance, the
knots $6_1$ and $9_{46}$ have the same Alexander polynomial, but
$\Augab(6_1,2)=2$ while $\Augab(9_{46},2)=5$.

More generally, knot contact homology seems to behave independently
of many known knot invariants. In particular, $HC_0$ can distinguish
knots which share many of the same ``classical'' invariants.

The knots $10_{25}$ and $\overline{10_{56}}$ (the mirror of
$10_{56}$, with the conventions of \cite{Rol}) are both alternating,
with the same HOMFLY polynomial and signature; hence they also share
the same Alexander polynomial, Jones polynomial, Khovanov invariant
\cite{Kho} (by a result of \cite{Lee}), and Ozsv\'ath--Szab\'o
invariant \cite{OSz}. On the other hand, we have
$\Augab(10_{25},7)=1$, while $\Augab(\overline{10_{56}},7)=2$.

Similarly, the knots $11_{255}^a$ and $11_{257}^a$, alternating
$11$-crossing knots from Conway's enumeration \cite{Con}, have the
same (two-variable) Kauffman polynomial, but
$\Augab(11_{255}^a,7)=1$ while $\Augab(11_{257}^a,7)=2$. We conclude
the following result.

\begin{proposition}
The invariant $HC_0$ for knots is not determined by any of the
following: Alexander polynomial, Jones polynomial, HOMFLY
polynomial, Kauffman polynomial, signature, Khovanov invariant, and
Ozsv\'ath--Szab\'o invariant.
\end{proposition}

%It is not clear whether $HC_0$ has the capability to distinguish
%mutants; for instance, it appears that $HC_0$ is the same for the
%Kinoshita--Terasaka knot and its Conway mutant.

%*********************************************************************
%*********************************************************************
\section{An alternate knot invariant}
\label{sec:alternate}

As explained in Section~\ref{sec:contact}, in the symplectic world,
the knot DGA counts holomorphic disks with certain boundary
conditions determined by the knot. This count is well defined modulo
$2$, and so we obtain a DGA over $\Z_2$. To lift this DGA to $\Z$
coefficients, however, we need a set of coherent orientations for
the relevant moduli spaces. The knot DGA we have been using
presumably corresponds to a choice of coherent orientations, but
there are others.

By experimenting with various choices of signs, one can derive
several inequivalent choices for DGAs over $\Z$ which are knot
invariants and agree with the original knot DGA modulo $2$, but only
two seem to give ``nice'' results for simple knots (eg, such that
$HC_0$ for small knots is of the form $\Z[x]/(p(x))$). One is the
knot DGA already described; the other we call the \textit{alternate
knot DGA}. This alternate theory lacks some of the nice properties
of the knot DGA---for instance, no obvious analogue of
Theorem~\ref{detthm} exists for the alternate DGA---but it seems to
give a sharper knot invariant than the usual DGA.

%*********************************************************************
\subsection{The alternate DGA}
\label{ssec:alternatedef}

We follow the definition of the knot DGA, making adjustments where necessary.

Let $\tilde{\phi}$ denote the homomorphism from $B_n$ to
$\Aut(\A_n)$ defined on generators by
\[
\varhom{\sigma_k}\co \left\{
\begin{array}{ccll}
a_{ki} & \mapsto & -a_{k+1,i} + a_{k+1,k}a_{ki}\phantom{999} & i\neq k,k+1 \\
a_{ik} & \mapsto & -a_{i,k+1} - a_{ik}a_{k,k+1} & i\neq k,k+1 \\
a_{k+1,i} & \mapsto & a_{ki} & i\neq k,k+1 \\
a_{i,k+1} & \mapsto & a_{ik} & i \neq k,k+1 \\
a_{k,k+1} & \mapsto & -a_{k+1,k} & \\
a_{k+1,k} & \mapsto & -a_{k,k+1} & \\
a_{ij} & \mapsto & a_{ij} & i,j \neq k,k+1.
\end{array}
\right.
\]
Note that $\tilde{\phi}$ is identical to $\phi$ modulo $2$, and that
$\phi$ and $\tilde{\phi}$ are conjugate through the automorphism of
$\A_n$ which sends $a_{ij}$ to $a_{ij}$ if $i<j$ and $-a_{ij}$ if
$i>j$.

Write $A=(a_{ij})$, $B=(b_{ij})$, $C=(c_{ij})$, $D=(d_{ij})$ as
usual, with the crucial difference that we set $a_{ii} = 0$, rather
than $a_{ii}=-2$, for all $i$. Denote by $\tilde{\phi}^{\ext}$ the
composition map $B_n \hookrightarrow B_{n+1}
\stackrel{\tilde{\phi}}{\rightarrow} \Aut(\A_n)$, and define
matrices $\varPhil{B}(A)$, $\varPhir{B}(A)$ by
\[
\varhom{B}^{\ext}(a_{i\dott}) = \sum_{j=1}^n (\varPhil{B}(A))_{ij} a_{j\dott}
\hspace{0.25in} \textrm{and} \hspace{0.25in}
 \varhom{B}^{\ext}(a_{\dott j}) =
\sum_{i=1}^n a_{\dott i} (\varPhir{B}(A))_{ij}.
\]
Then the \textit{alternate knot DGA} of $B \in B_n$ is generated, as
before, by $\{a_{ij}\,|\,1\leq i\neq j\leq n\}$ in degree $0$,
$\{b_{ij},c_{ij}\,|\,1\leq i,j\leq n\}$ in degree $1$, and
$\{d_{ij}\,|\,1\leq i,j\leq n\}$ and $\{e_i\,|\,1\leq i\leq n\}$ in
degree $2$, with
\begin{align*}
\d A &= 0 \\
\d B &= (1 + \varPhil{B}(A)) \cdot A \\
\d C &= A \cdot (1 + \varPhir{B}(A)) \\
\d D &= B \cdot (1 + \varPhir{B}(A)) - (1 + \varPhil{B}(A)) \cdot C \\
\d e_i &= (B - \varPhil{B}(A) \cdot C)_{ii}.
\end{align*}
Note the difference in signs from Definition~\ref{knotdga}. Also,
let the \textit{alternate contact homology} of a knot $K$, written
$\varHC_*(K)$, be the homology of the alternate knot DGA of any
braid closing to $K$. As mentioned before, the alternate knot DGA
and contact homology are identical modulo $2$ to the usual knot DGA
and contact homology.

\begin{theorem}
The stable tame isomorphism class of the alternate knot DGA, as well as $\varHC_*$, are
\label{varknotthm}
knot invariants.
\end{theorem}

\noindent The proof of Theorem~\ref{varknotthm} is entirely
analogous to the proofs from Sections~\ref{sec:HC}
and~\ref{sec:fullproof}, and is omitted here.

Many of the results from Section~\ref{sec:properties} still hold in
this context. We can define conjugation on the alternate knot DGA by
$\overline{a_{ij}} = -a_{ji}$, $\overline{b_{ij}} = -c_{ji}$,
$\overline{c_{ij}} = -b_{ji}$, $\overline{d_{ij}} = d_{ji}$, and
$\overline{e_i} = d_{ii}-e_i$; then the differential $\d$ in the
alternate knot DGA commutes with conjugation. Hence we can form the
\textit{alternate abelian knot DGA} similarly to
Definition~\ref{abeliandef}, but with $a_{ji}=-a_{ij}$ and
$c_{ij}=-b_{ji}$. None of the alternate invariants distinguish
between knot mirrors or inverses.

%*********************************************************************
\subsection{Computations of $\varHC_0(K)$ for small knots}
\label{ssec:varsmallknots}

One interesting question is whether an analogue of
Theorem~\ref{detthm} holds for $\varHC_0$. In all computations
performed by the author, there is a surjection from $\varHC_0(K)$ to
some ring of the form $\Z[x]/(p(x))$, where $p$ is monic of degree
$\frac{n(K)+1}{2}$ and $p(x) \equiv p_{(n(K)+1)/2}(x) \pmod 2$. Why
this should be true in general is unclear.

Nevertheless, for small knots $K$, we have $\varHC_0(K) \cong
\Z[x]/(p(x))$ for some polynomial $p$ with $\deg(p) =
(|\Delta_K(-1)|+1)/2$. In contrast to the case of $HC_0$, however,
the polynomial $p$ depends on more than just $|\Delta_K(-1)|$. For
example, the knots $4_1$ and $5_1$ both satisfy $|\Delta_K(-1)| = 5$
and $HC_0(K) \cong \Z[x]/(x^3-x^2-3x+2)$, but $\varHC_0(4_1) \cong
\Z[x]/(x^3-x^2+x)$, while $\varHC_0(5_1) \cong \Z[x]/(x^3-x^2-x)$.
See Table~\ref{knottable} for more examples.

\begin{table}[ht!]\small
\[
\begin{array}{|l||l|r|}
\hline \heighten K & \hfill p(x)~\textrm{with}~\varHC_0(K)\cong \Z[x]/(p(x)) \hfill &
\hfill \disc(p) \hfill \rule{0pt}{12pt} %& \deg(p)& |\Delta_K(-1)|
\\
\hline \hline 0_1 & x & 1 \heighten %& 1 & 1
\\
\hline 3_1 & x^2-x & 1 \heighten %& 2 & 3
\\
\hline 4_1 & x^3-x^2+x & -3 \heighten %& 3 & 5
\\
\hline 5_1 & x^3-x^2-x & 5 \heighten %& 3 & 5
\\
\hline 5_2 & x^4-x^3+x & -23 \heighten %& 4 & 7
\\
\hline 6_1 & x^5-x^4+x^3+x & 257 \heighten %& 5 & 9
\\
\hline 6_2 & x^6-x^5-2x^4+x^3+x^2+x & 1777 \heighten %& 6 & 11
\\
\hline 6_3 & x^7-x^6-x^5+2x^4-x^2+x & -10571 \heighten %& 7 & 13
\\
\hline 7_1 & x^4-x^3-2x^2+x & 49 \heighten %& 4 & 7
\\
\hline 7_2 & x^6-x^5+x^3+x^2-x & 4409 \heighten %& 6 & 11
\\
\hline 7_3 & x^7-x^6-3x^5+2x^4+2x^3+x^2-x & 78301 \heighten %& 7 & 13
\\
\hline 7_4 & x^8-x^7+4x^6-3x^5+4x^4-2x^3+x & -1166607 \heighten %& 8 & 15
\\
\hline 7_5 & x^9-x^8-x^7+2x^6+x^5-2x^4+2x^2-x & -4690927 \heighten %& 9 & 17
\\
\hline 7_6 & x^{10}-x^9+2x^8-x^7+3x^6-x^5+2x^4+x^2+x & 90320393 \heighten %& 10 & 19
\\
\hline 7_7 & x^{11}-x^{10}+3x^9-2x^8+4x^7-3x^6+3x^5-2x^4+x^3-x^2+x
& -932501627 \heighten %& 11 & 21
\\
\hline
\end{array}
\]
\caption{The invariant $\varHC_0(K)$ for prime knots with seven or fewer crossings,
and the discriminant $\disc(p)$. Knot notation is as in \cite{Rol}.
}
\label{knottable}
\end{table}

\begin{proposition}
If $K$ is a prime knot with seven or fewer crossings, then we can write
$\varHC_0(K) \cong \Z[x]/(p(x))$ with $\deg(p) = (|\Delta_K(-1)|+1)/2$.
\label{varsmallknots}
Furthermore, $p$ can be expressed as
$p=\tilde{p}_{(|\Delta_K(-1)|+1)/2}$ for some sequence
$\{\tilde{p}_m\}$ with $\tilde{p}_0(x)=\tilde{p}_1(x)=x$ and
$\tilde{p}_{m+1}(x) = x\tilde{p}_m(x) \pm \tilde{p}_{m-1}(x)$ for
each $m$; in particular, $p \equiv p_{(|\Delta_K(-1)|+1)/2} \pmod
2$, with $\{p_m\}$ as in Section~\ref{ssec:det}.
\end{proposition}

\noindent Proposition~\ref{varsmallknots} can be verified by direct
computation using the author's program in \textit{Mathematica}; see
Table~\ref{knottable}.

As with $HC_0$, it is easy to compute $\varHC_0$ explicitly for any
torus knot $T(2,n)$.

\begin{proposition}
For $k \geq 1$, $HC_0(T(2,2k-1)) \cong \Z[x]/(\tilde{p}_k(x))$,
where \label{var31}
\[
\tilde{p}_k(x) = \sum_{j=0}^{k-1} (-1)^{\lceil j/2 \rceil}
\binom{k-1-\lceil j/2 \rceil}{\lfloor j/2 \rfloor} x^{k-j}.
\]
\end{proposition}

\noindent Note that the sequence $\{\tilde{p}_k\}$ in
Proposition~\ref{var31} satisfies $\tilde{p}_0(x)=\tilde{p}_1(x)=x$
and $\tilde{p}_{m+1}(x) = x\tilde{p}_m(x) - \tilde{p}_{m-1}(x)$, and
that $\tilde{p}_m(x) = x p_m(x)/(x+2)$ for all $m$.

In general, when $\varHC_0(K)$ is of the form $\Z[x]/(p(x))$ with
$p$ monic, the polynomial $p(x)$ is not unique; however, we can
deduce two useful invariants of the isomorphism class of
$\varHC_0(K)$ from $p(x)$. The first is the degree of $p$, which is
simply the vector space dimension of $\varHC_0(K;\Q)$. The second is
the discriminant $\disc(p)$, the well-known integer invariant which
can be defined as $\prod_{i<j} (r_i-r_j)^2$, where $\{r_i\}$ is the
set of roots of $p$.

\begin{table}[ht!]\small
\[
\begin{array}{|l||l|l|}
\hline K &
\hfill \heighten p(x)~\textrm{with}~HC_0(K)\cong \Z[x]/(p(x)) \hfill &
\hfill p(x)~\textrm{with}~\varHC_0(K)\cong \Z[x]/(p(x)) \hfill \rule{0pt}{12pt}
\\
\hline \hline 0_1 & x-2 & x \heighten %& 1 & 1
\\
\hline 3_1 & x^2-x-2 & x^2-x \heighten %& 2 & 3
\\
\hline 4_1 & x^3-x^2-3x+2 & x^3-x^2+x \heighten %& 3 & 5
\\
\hline 5_1 & x^3-x^2-3x+2 & x^3-x^2-x \heighten %& 3 & 5
\\
\hline 5_2 & x^4-x^3-4x^2+3x+2 & x^4-x^3+x \heighten %& 4 & 7
\\
\hline 6_1 & x^5-x^4-5x^3+4x^2+5x-2 & x^5-x^4+x^3+x \heighten %& 5 & 9
\\
\hline 6_2 & x^6-x^5-6x^4+5x^3+9x^2-5x-2 & x^6-x^5-2x^4+x^3+x^2+x \heighten %& 6 & 11
\\
\hline 6_3 & x^7-x^6-7x^5+6x^4+14x^3-9x^2-7x+2 & x^7-x^6-x^5+2x^4-x^2+x \heighten
\\
\hline 7_1 & x^4-x^3-4x^2+3x+2 & x^4-x^3-2x^2+x \heighten %& 4 & 7
\\
\hline 7_2 & x^6-x^5-6x^4+5x^3+9x^2-5x-2 & x^6-x^5+x^3+x^2-x \heighten %& 6 & 11
\\
\hline 7_3 & x^7-x^6-7x^5+6x^4+14x^3-9x^2-7x+2 & x^7-x^6-3x^5+2x^4+2x^3+x^2-x \heighten %& 7 & 13
\\
\hline
\end{array}
\]
\caption{Comparison between $HC_0$ and $\varHC_0$ for some small knots.
}
\label{knottable2}
\end{table}

Table~\ref{knottable} gives $\varHC_0(K)$ for all prime knots $K$
with seven or fewer crossings, while Table~\ref{knottable2}
contrasts the invariants $HC_0$ and $\varHC_0$. Note that the
discriminant on $\varHC_0$ allows us to distinguish between knots
which have the same $HC_0$; in other words, lifting
$\varHC_0(K;\Z_2)=HC_0(K;\Z_2)$ to $\varHC_0(K)$ does give new
information. This would seem to provide an instance where
considering orientations in contact homology yields results not
given by the unoriented theory.

%*********************************************************************
%*********************************************************************

\end{document}